\documentclass[11pt]{article}
\usepackage{latexsym}

\newcommand{\eproof}{\mbox{\ }\hfill $\Box$ \par \vskip 10pt}

\newtheorem{Theorem}{Theorem}[section]
\newtheorem{lemma}[Theorem]{Lemma}
\newtheorem{prop}[Theorem]{Proposition}

\begin{document}

\title{Semi-classical dispersive estimates}

\author{{\sc Fernando Cardoso, Claudio Cuevas and Georgi Vodev\thanks{Corresponding author}}}

\date{}

\maketitle

\noindent
{\bf Abstract.} We prove dispersive estimates for the wave group $e^{it\sqrt{P(h)}}$ and the Schr\"odinger group $e^{itP(h)}$, where $P(h)$ is a
self-adjoint, elliptic second-order differential operator depending on a parameter $0<h\le 1$, which is supposed to be a short-range perturbation of $-h^2\Delta$, $\Delta$ being the Euclidean Laplacian. In particular, applications are made to non-trapping metric perturbations and to perturbations by a magnetic potential. 

\setcounter{section}{0}
\section{Introduction and statement of results}

Denote by $P_0(h)$ the self-adjoint realization of $-h^2\Delta$ on $L^2({\bf R}^n)$, $n\ge 2$, and let 
$\varphi\in C_0^\infty((0,+\infty))$ be independent of $h$. It is well known (see the appendix) that the free wave and 
Schr\"odinger groups satisfy the following dispersive estimates
$$\left\|\langle x\rangle^{-\sigma}e^{it\sqrt{P_0(h)}}\varphi(P_0(h))\langle x\rangle^{-\sigma}\right\|_{L^1\to L^\infty}\le Ch^{-n-\sigma}t^{-\frac{n-1}{2}-\sigma},
\eqno{(1.1)}$$
$$\left\|\langle x\rangle^{-\sigma}e^{itP_0(h)}\varphi(P_0(h))\langle x\rangle^{-\sigma}\right\|_{L^1\to L^\infty}\le Ch^{-n-\sigma}t^{-\frac{n}{2}-\sigma},
\eqno{(1.2)}$$
for all $t>0$, $0<h\le 1$, $\sigma\ge 0$, 
with a constant $C>0$ independent of $t$ and $h$. The purpose of this paper is to prove analogues of (1.1) and 
(1.2) for more general second-order operators of the form
$$P(h)=\sum_{i,j=1}^n{\cal D}_{x_i}a_{ij}(x,h){\cal D}_{x_j}+\sum_{j=1}^n\left(b_j(x,h){\cal D}_{x_j}+
{\cal D}_{x_j}b_j(x,h)\right)+V(x,h),$$
with real-valued coefficients $a_{ij},b_j\in C^1({\bf R}^n)$ and $V\in L^\infty({\bf R}^n)$, 
where ${\cal D}_{x_j}:=-ih\partial_{x_j}$, $0<h\le 1$ is a semi-classical parameter (not necessarilly small). 
More precisely, the coefficients are of the form
$$a_{ij}(x,h)=a_{ij}^0(x)+ha_{ij}^1(x,h),\,b_j(x,h)=b_j^0(x)+hb_j^1(x,h),\,V(x,h)=V^0(x)+hV^1(x,h),$$
where $a_{ij}^0,b_j^0,V^0\in C^1({\bf R}^n)$ are independent of $h$, and $a_{ij}^1,b_j^1,V^1\in L^\infty({\bf R}^n)$ 
uniformly in $h$. So, the principal symbol of $P(h)$ is given by
$$p(x,\xi)=\sum_{i,j=1}^na_{ij}^0(x)\xi_i\xi_j+2\sum_{j=1}^nb_j^0(x)\xi_j+V^0(x).$$
We suppose that this operator admits a self-adjoint realization on the Hilbert space $L^2({\bf R}^n)$ 
(which will be again denoted by $P(h)$) satisfying the ellipticity condition
$$\sum_{0\le|\alpha|\le 2}\left\|{\cal D}_x^\alpha(P(h)\pm i)^{-1}\right\|_{L^2\to L^2}\le C,\eqno{(1.3)}$$
with a constant $C>0$ independent of $h$. We also suppose that $P(h)$ is a short-range perturbation of $P_0(h)$, namelly
$$\sum_{i,j=1}^n\left|\partial_x^\alpha\left(a_{ij}(x,h)-a^\flat_{ij}\right)\right|+
\sum_{j=1}^n\left|\partial_x^\alpha b_j(x,h)\right|+\left|\partial_x^\alpha V^0(x)\right|+
\left|V^1(x,h)\right|\le C\langle x\rangle^{-\delta},\eqno{(1.4)}$$
for $0\le|\alpha|\le 1$, 
with constants $C>0$, $\delta>1$ independent of $h$, where $a^\flat_{ij}=1$ if $i=j$, $a^\flat_{ij}=0$ if $i\neq j$. 
When $n\ge 4$ we suppose that there exists a sufficiently small constant $\gamma>0$, independent of $h$, such that
$$\sup_{x\in{\bf R}^n}\sum_{i,j=1}^n\left|a_{ij}(x,h)-a^\flat_{ij}\right|\le\gamma.\eqno{(1.5)}$$
Given a $z\in{\bf C}$, Im$\,z\neq 0$, set
$$R_s(z,h):=\langle x\rangle^{-s}\left(P(h)-z\right)^{-1}\langle x\rangle^{-s}.$$
Finally, we suppose that there exist an energy level $E>0$ and a constant $0<\varepsilon_0<E$, both independent of $h$, 
so that for every $z\in [E-\varepsilon_0,E+\varepsilon_0]$, $s>1/2$, the limits 
$$R_s^{\pm}(z,h):=\lim_{\varepsilon\to 0^+}R_s(z\pm i\varepsilon,h):L^2\to L^2$$
exist as continuous functions in $z$ and satisfy the bound
$$\left\|R_s^{\pm}(z,h)\right\|_{L^2\to L^2}\le C\mu(h),\quad \forall z\in[E-\varepsilon_0,E+\varepsilon_0],\eqno{(1.6)}$$
with a constant $C>0$ and a function $\mu(h)\ge h^{-1}\ge 1$. If the coefficients are smooth and if $E$ 
is a non-trapping energy level, i.e. all
bicharacteristics belonging to $\{(x,\xi)\in T^*{\bf R}^n:\,p(x,\xi)=E\}$ escape to infinity, 
it is well known that (1.6) holds with $\mu(h)=h^{-1}$
provided $\varepsilon_0$ is taken small enough independent of $h$. More generally, it is proved in 
\cite{kn:NZ} that (1.6) holds with $\mu(h)=h^{-1}\log\left(h^{-1}\right)$ if all periodic bicharacteristics 
belonging to $\{(x,\xi)\in T^*{\bf R}^n:\,p(x,\xi)=E\}$ are of hyperbolic type satisfying a topological condition. 
On the other hand, without any geometrical condition we have that $\mu(h)=e^{\beta/h}$, $\beta>0$ a constant, 
still for smooth coefficients (e.g. see \cite{kn:Bu}, \cite{kn:CV1}). Hence, in this case the function $\mu$
satisfies
$$h^{-1}\le\mu(h)\le e^{\beta/h},\quad \beta>0.\eqno{(1.7)}$$
It is largely expected that (1.7) holds true under the assumptions above.

Let $\varphi\in C_0^\infty((E-\varepsilon_0,E+\varepsilon_0))$ be independent of $h$. 
In the present paper we are interested in bounding from above uniformly in $h$ the following quantities
$$A_1(h,\sigma)=h^{n+\sigma}\sup_{f\in \langle x\rangle^{-s}L^2,\,\|\langle x\rangle^{s}f\|_{L^2}=1}\,\sup_{t>0}\,t^{(n-1)/2+\sigma}
\left\|\langle x\rangle^{-\sigma}e^{it\sqrt{P(h)}}\varphi(P(h))\langle x\rangle^{-\sigma}f\right\|_{L^\infty},$$
where $s>n/2$, 
$$A_2(h,\sigma)=h^{n+\sigma}\sup_{f\in L^1,\,\|f\|_{L^1}=1}\,\sup_{t>0}\,t^{(n-1)/2+\sigma}\left\|\langle x\rangle^{-\sigma}e^{it\sqrt{P(h)}}
\varphi(P(h))\langle x\rangle^{-\sigma}f\right\|_{L^\infty},$$
 $$B_1(h,\sigma)=h^{n+\sigma}\sup_{f\in \langle x\rangle^{-s}L^2,\,\|\langle x\rangle^{s}f\|_{L^2}=1}\,\sup_{t>0}\,t^{n/2+\sigma}
\left\|\langle x\rangle^{-\sigma}e^{itP(h)}\varphi(P(h))\langle x\rangle^{-\sigma}f\right\|_{L^\infty},$$
where $s>(n+1)/2$,
 $$B_2(h,\sigma)=h^{n+\sigma}\sup_{f\in L^1,\,\|f\|_{L^1}=1}\,\sup_{t>0}\,t^{n/2+\sigma}\left\|\langle x\rangle^{-\sigma}e^{itP(h)}\varphi(P(h))\langle x\rangle^{-\sigma}f\right\|_{L^\infty}.$$
In view of the estimates (1.1) and (1.2), we have that in the case of the free operator $P_0(h)$ all these 
quantities are bounded by a constant independent of $h$. In the present paper we will show that in the general 
case of the operator $P(h)$ these quantities can be bounded from above 
in terms of the function $\mu(h)$, provided $\delta$ is big enough. To state our main result more precisely, 
we define the number $\nu\in\{0,1,2\}$ as follows: $\nu =2$ if $a^0_{ij}(x)-a^\flat_{ij}\equiv b_j^0(x)
\equiv V^0(x)\equiv 0$, and the functions $a^1_{ij}(x,h)$, 
$b_j^1(x,h)$, $V^1(x,h)$ are $O(h)$ as $h\to 0$; $\nu =1$ if $a^0_{ij}(x)-a^\flat_{ij}\equiv b_j^0(x)\equiv 
V^0(x)\equiv 0$; $\nu=0$ otherwise.
In other words, the quantity $2-\nu$ can be viewed as the order of the perturbation $P(h)-P_0(h)$. We have the following

\begin{Theorem} Suppose the conditions (1.3)-(1.6) satisfied with $\delta>\frac{n+2}{2}+\sigma$ in the case of the wave group and 
$\delta>\frac{n+3}{2}+\sigma$ in the case of the Schr\"odinger group with some $\sigma\ge 0$. Then the following bounds hold true:
$$A_1(h,\sigma)\le C_\varepsilon h^{\nu+\sigma+\frac{n-1}{2}}\mu(h)^{\frac{n+1}{2}+\sigma+\varepsilon}+C,\eqno{(1.8)}$$
$$A_2(h,\sigma)\le C_\varepsilon h^{2\nu+\sigma -\frac{3}{2}}\mu(h)^{\frac{n+1}{2}+\sigma+\varepsilon}+Ch^{\nu -\frac{n+1}{2}}+C,\eqno{(1.9)}$$
$$B_1(h,\sigma)\le C_\varepsilon h^{\nu+\sigma+\frac{n-1}{2}}\mu(h)^{\frac{n+2}{2}+\sigma+\varepsilon}+C,\eqno{(1.10)}$$
$$B_2(h,\sigma)\le C_\varepsilon h^{2\nu+\sigma -\frac{3}{2}}\mu(h)^{\frac{n+2}{2}+\sigma+\varepsilon}+Ch^{\nu -\frac{n+1}{2}}+C,\eqno{(1.11)}$$
for every $0<\varepsilon\ll 1$.
\end{Theorem}

We will apply these estimates to operators of the form $P(h)=h^2G$, where $G$ is the self-adjoint realization of a 
second-order operator of the form
$$G=-\sum_{i,j=1}^n\partial_{x_i}a_{ij}(x)\partial_{x_j}+i\sum_{j=1}^n\left(b_j(x)\partial_{x_j}+\partial_{x_j}b_j(x)\right)+V(x),$$
with real-valued coefficients $a_{ij},b_j\in C^1({\bf R}^n)$, $V\in L^\infty({\bf R}^n)$ independent of $h$, satisfying
$$\sum_{i,j=1}^n\left|\partial_x^\alpha\left(a_{ij}(x)-a^\flat_{ij}\right)\right|+\sum_{j=1}^n\left|\partial_x^\alpha b_j(x)
\right|+\left|V(x)\right|\le C\langle x\rangle^{-\delta},\eqno{(1.12)}$$
for $0\le|\alpha|\le 1$, 
with constants $C>0$, $\delta>1$. In other words, $G$ is supposed to be a short-range perturbation of the 
self-adjoint realization, $G_0$, of
the free Laplacian $-\Delta$. When $n\ge 4$ we suppose that there exists a sufficiently small constant $\gamma>0$ such that
$$\sup_{x\in{\bf R}^n}\sum_{i,j=1}^n\left|a_{ij}(x)-a^\flat_{ij}\right|\le\gamma.\eqno{(1.13)}$$
We also suppose that $G$ is elliptic, that is,
$$\partial_x^\alpha (G\pm i)^{-1}\in {\cal L}(L^2),\eqno{(1.14)}$$
for all $|\alpha|\le 2$, where ${\cal L}(L^2)$ denotes the set of the bounded operators on $L^2$.  
We finally suppose that there exist constants $C,\lambda_0>0$ and $k\ge 0$ such that
$$\left\|\langle x\rangle^{-s}\left(G-\lambda^2\pm i0\right)^{-1}\langle x\rangle^{-s}\right\|_{L^2\to L^2}\le C\lambda^{-1+k},
\quad \forall \lambda\ge\lambda_0,\,s>1/2.\eqno{(1.15)}$$
This implies that the operator $P(h)=h^2G$ satisfies (1.6) with $\mu(h)=h^{-1-k}$. Set
$$p_n(\sigma)=\max\left\{0,\frac{n+1}{2}-\nu,\frac{n+4}{2}-2\nu+\frac{k(n+1)}{2}+k\sigma\right\},$$
$$q_n(\sigma)=\max\left\{0,\frac{n+5}{2}-2\nu+\frac{k(n+2)}{2}+k\sigma\right\},$$
where $\sigma\ge 0$ and $2-\nu$ is the order of the differential operator $G-G_0$. 
Let $\chi\in C^\infty({\bf R})$, supp$\,\chi\subset(\lambda_0,+\infty)$, $\chi(\lambda)=1$ for $\lambda\ge \lambda_0+1$. 
As a consequence of the above theorem we get the following $\langle x\rangle^{\sigma}L^1\to \langle x\rangle^{-\sigma}L^\infty$ dispersive estimates for the perturbed wave 
(resp. Schr\"odinger) group with a loss of $p_n(\sigma)$ (resp. $q_n(\sigma)+\varepsilon$) derivatives.

\begin{Theorem} Suppose the conditions (1.12)-(1.15) satisfied with $\delta>\frac{n+2}{2}+\sigma$ in the case of 
the wave group and $\delta>\frac{n+3}{2}+\sigma$ in the case of the Schr\"odinger group with some $\sigma\ge 0$. Then, the following dispersive 
estimates hold true:
$$\left\|\langle x\rangle^{-\sigma}e^{it\sqrt{G}}(\sqrt{G})^{-\frac{(n+1)}{2}-p_n(\sigma)-\varepsilon}\chi(\sqrt{G})\langle x\rangle^{-\sigma}\right\|_{L^1\to L^\infty}
\le C_\varepsilon|t|^{-\frac{n-1}{2}-\sigma},\quad\forall t\neq 0,\eqno{(1.16)}$$
$$\left\|\langle x\rangle^{-\sigma}e^{itG}(\sqrt{G})^{\sigma-q_n(\sigma)-\varepsilon}\chi(\sqrt{G})\langle x\rangle^{-\sigma}\right\|_{L^1\to L^\infty}\le C_\varepsilon|t|^{-\frac{n}{2}-\sigma},
\quad\forall t\neq 0,\eqno{(1.17)}$$
for every $0<\varepsilon\ll 1$. Moreover, if $k<1$ and
$$\delta>\frac{n+3}{2}+\frac{q_n(0)}{1-k},$$
then for all $\sigma$ satisfying
$$\frac{q_n(0)}{1-k}<\sigma<\delta-\frac{n+3}{2},$$
we have the estimate
$$\left\|\langle x\rangle^{-\sigma}e^{itG}\chi(\sqrt{G})\langle x\rangle^{-\sigma}\right\|_{L^1\to L^\infty}\le C|t|^{-\frac{n}{2}-\sigma},
\quad\forall t\neq 0.\eqno{(1.18)}$$
\end{Theorem}

In the particular case of non-trivial non-trapping metric perturbations we have (1.15) with $k=0$ as well as $\nu=0$, 
so $p_n(\sigma)=\frac{n+4}{2}$, 
$q_n(\sigma)=\frac{n+5}{2}$. Thus, in this case we obtain $\langle x\rangle^{\sigma}L^1\to \langle x\rangle^{-\sigma}L^\infty$ dispersive estimates for the perturbed wave 
(resp. Schr\"odinger) group with a loss of $\frac{n+4}{2}$ (resp. $\frac{n+5}{2}+\varepsilon$) derivatives. 
The same conclusion remains true for more general metric perturbations with infinitely many periodic geodesics of hyperbolic type. 
Indeed, for such perturbations the bound (1.15) with $k=\varepsilon$, $\forall 0<\varepsilon\ll 1$, has been proved in 
\cite{kn:NZ} under some natural topological conditions. We get a better result for perturbations by a
magnetic potential, namely for operators of the form
$$G=\left(i\nabla+b(x)\right)^2+V(x),$$
where $b(x)=\left(b_1(x),...,b_n(x)\right)\in C^1\left({\bf R}^n;{\bf R}^n\right)$ is a vector-valued function and 
$V\in L^\infty\left({\bf R}^n;{\bf R}\right)$. When $n\ge 3$ it is proved in \cite{kn:EGS} (see Proposition 4.3) 
that in this case (1.15) holds with $k=0$. Since $\nu=1$, we have in this case $p_n(\sigma)=\frac{n}{2}$, $q_n(\sigma)=\frac{n+1}{2}$. 
When $b(x)\equiv 0$, we have $\nu=2$ and hence in this case $p_n(\sigma)=q_n(\sigma)=\frac{n-3}{2}$.
This latter case, however, has already been studied in \cite{kn:CV2}, \cite{kn:CCV1}, \cite{kn:V1}, \cite{kn:V2} 
under a little bit weaker assumption on the potential $V$. 

To our best knowledge, it is the first time dispersive estimates are proved for perturbations different from a potential. 
Our estimates are not
optimal (i.e. we are obliged to loose derivatives), but one could hardly do better without assuming a stronger 
regularity of the coefficients. Indeed, it was shown in \cite{kn:GV} in the context of the Schr\"odinger equation 
with a potential that if $n\ge 4$, it is not possible to have optimal $L^1\to L^\infty$ dispersive estimates for potentials 
$V\in C_0^k({\bf R}^n)$, $\forall k<\frac{n-3}{2}$. In contrast, when $n\le 3$ no
regularity of the potential is needed in order to have optimal $L^1\to L^\infty$ dispersive estimates for both the wave and 
the Schr\"odinger groups (e.g. see \cite{kn:M} when $n=2$ and \cite{kn:DP}, \cite{kn:G} when $n=3$). When $n\ge 4$ 
it is expected that optimal dispersive estimates hold true for potentials $V\in C^{\frac{n-3}{2}}({\bf R}^n)$. 
Indeed, such results have been recently proved in \cite{kn:EG} when $n=5,7$, (see also \cite{kn:CCV}) in the
case of the Schr\"odinger equation and in \cite{kn:CV3} when $4\le n\le 7$ in the case of the wave equation. 
For potentials with stronger regularity optimal dispersive estimates were proved in \cite{kn:B} in the case of 
the wave equation with Schwartz class potentials and in \cite{kn:JSS} (see also \cite{kn:MV}) in the case of the 
Schr\"odinger equation with potentials satisfying $\widehat V\in L^1$. To our best knowledge, no optimal dispersive 
estimates have been proved so far in the more general context of the operator $G$ above when the function
$$\sum_{i,j=1}^n\left|a_{ij}(x)-a^\flat_{ij}\right|+\sum_{j=1}^n\left|b_j(x)\right|$$
is not identically zero, even if we suppose that $a_{ij}-a^\flat_{ij},b_j,V\in C_0^\infty({\bf R}^n)$. In general, 
proving optimal dispersive estimates turns out to be a very tough problem.

To prove the main result we extend to more general perturbations the method developed in \cite{kn:V1}, \cite{kn:V2}, 
which consists of deriving the dispersive estimates from decay estimates on weighted $L^2$ spaces. This analysis is based 
on a careful study of the regularity of the resolvent on weighted $L^2$ spaces (see Proposition 3.2 below). 
Note that the assumption (1.5) is only used in the proof of the estimate (2.5) which plays a crucial role in our approach. 
It might be possible, however, that (2.5) could hold without (1.5). It becomes clear from the proof
that the reason why we need (1.5) is due to the fact that when $n\ge 4$ the singularity at zero of the Hankel 
functions is too strong, which in turn implies a very strong singularity on the diagonal of the kernel of the free resolvent. 
Consequently, $(P_0(h)-z)^{-1}:L^2\to L^\infty$, ${\rm Im}\,z\neq 0$, is no longer bounded when $n\ge 4$. 
This difficulty is overcome by Lemma 2.2 below. Note finally that we expect that the above estimates hold true for
$\delta>\frac{n+1}{2}$, but this is much harder to prove especially in the case of the Schr\"odinger group. 

\section{Study of the operator $\varphi(P(h))$}

In this section we will prove the following

\begin{prop} Assume (1.3),(1.4) and (1.5) fulfilled. Then, for all $0\le s,s_1,s_2\le\delta$, $s_1+s_2\le\delta$, we have the bounds
$$\left\|\langle x\rangle^{-s}\varphi(P_0(h))\langle x\rangle^s\right\|_{L^2\to L^2}\le C,\eqno{(2.1)}$$
$$\left\|\langle x\rangle^{-s}\varphi(P(h))\langle x\rangle^s\right\|_{L^2\to L^2}\le C,\eqno{(2.2)}$$
$$\left\|\langle x\rangle^{s_1}\left(\varphi(P(h))-\varphi(P_0(h))\right)\langle x\rangle^{s_2}\right\|_{L^2\to L^2}
\le Ch^\nu,\eqno{(2.3)}$$ 
$$\left\|\langle x\rangle^{s_1}\varphi(P_0(h))\left(P(h)-P_0(h)\right)\varphi(P(h))\langle x\rangle^{s_2}
\right\|_{L^2\to L^2}\le Ch^{\nu},\eqno{(2.4)}$$
$$\left\|\left(\varphi(P(h))-\varphi(P_0(h))\right)\langle x\rangle^s\right\|_{L^2\to L^\infty}\le Ch^{\nu-n/2},\eqno{(2.5)}$$
with a constant $C>0$ independent of $h$.
\end{prop}

{\it Proof.} The estimate (2.1) is well known, while (2.2) follows from (2.1) and (2.3). It is also easy to see that (2.4) 
follows from
(2.3). To prove (2.3) and (2.5) we will use the Helffer-Sj\"ostrand formula
$$\varphi(P(h))=\frac{1}{\pi}\int_{\bf C}\frac{\partial\widetilde\varphi}{\partial\bar z}(z)\left(P(h)-z\right)^{-1}L(dz),
\eqno{(2.6)}$$
where $L(dz)$ denotes the Lebesgue measure on ${\bf C}$, $\widetilde\varphi\in C_0^\infty({\bf C})$ is an almost 
analytic continuation
of $\varphi$ supported in a small complex neighbourhood of supp$\,\varphi$ and satisfying
$$\left|\frac{\partial\widetilde\varphi}{\partial\bar z}(z)\right|\le C_N|{\rm Im}\,z|^N,\quad\forall N\ge 1.\eqno{(2.7)}$$
It is well known that the free resolvent satisfies the estimate (e.g. see the proof of Lemma 2.3 of \cite{kn:V1})
$$\left\|\langle x\rangle^{-s}{\cal D}_x^\alpha\left(P_0(h)-z\right)^{-1}\langle x\rangle^s\right\|_{L^2\to L^2}
\le C\left|{\rm Im}\,z\right|^{-q},\quad |\alpha|\le 2,\eqno{(2.8)}$$
for $z\in{\rm supp}\,\widetilde\varphi$, ${\rm Im}\,z\neq 0$, with a constant $C>0$ independent of $z$ and $h$. 
Let us see that a similar estimate
holds true for the perturbed resolvent. Recall first that by assumption
$$P(h)-P_0(h)=h^\nu\sum_{|\alpha|\le 2}r_\alpha(x,h){\cal D}_x^\alpha,\eqno{(2.9)}$$
with coefficients satisfying
$$\left|r_\alpha(x,h)\right|\le C\langle x\rangle^{-\delta},\eqno{(2.10)}$$
with a constant $C>0$ independent of $x$ and $h$. Note also that (1.3) implies
$$\left\|{\cal D}_x^\alpha\left(P(h)-z\right)^{-1}\right\|_{L^2\to L^2}\le C\left|{\rm Im}\,z\right|^{-1},
\quad |\alpha|\le 2,\eqno{(2.11)}$$
for $z\in{\rm supp}\,\widetilde\varphi$, ${\rm Im}\,z\neq 0$. Using (2.8)-(2.11) together with the resolvent identity, we obtain
$$\left\|\langle x\rangle^{-s}{\cal D}_x^\alpha\left(P(h)-z\right)^{-1}\langle x\rangle^s\right\|_{L^2\to L^2}
\le \left\|\langle x\rangle^{-s}{\cal D}_x^\alpha\left(P_0(h)-z\right)^{-1}\langle x\rangle^s\right\|_{L^2\to L^2}$$
 $$+C\sum_{|\beta|\le 2}\left\|{\cal D}_x^\alpha\left(P(h)-z\right)^{-1}\right\|_{L^2\to L^2}\left\|
\langle x\rangle^{-\delta}{\cal D}_x^\beta\left(P_0(h)-z\right)^{-1}\langle x\rangle^s
\right\|_{L^2\to L^2}$$ $$\le C\left|{\rm Im}\,z\right|^{-q-1}.\eqno{(2.12)}$$
 On the other hand, using (2.6), (2.9), (2.10) and the resolvent identity, we get
 $$\left\|\langle x\rangle^{s_1}\left(\varphi(P(h))-\varphi(P_0(h))\right)\langle x\rangle^{s_2}
\right\|_{L^2\to L^2}$$ $$\le Ch^{\nu}\sum_{|\alpha|\le 2}
\int_{\bf C}\left|\frac{\partial\widetilde\varphi}{\partial\bar z}(z)\right|\left\|\langle x\rangle^{s_1}
\left(P(h)-z\right)^{-1}\langle x\rangle^{-s_1}\right\|_{L^2\to L^2}$$ $$\times\left\|\langle x
\rangle^{s_1-\delta}{\cal D}_x^\alpha\left(P_0(h)-z\right)^{-1}\langle x\rangle^{s_2}\right\|_{L^2\to L^2}L(dz).\eqno{(2.13)}$$
Clearly, (2.3) follows from (2.7), (2.8), (2.12) and (2.13).

To prove (2.5) we will first consider the case $n=2,3$. Then it is well known that the free resolvent satisfies the estimate
$$\left\|\left(P_0(h)-z\right)^{-1}\right\|_{L^2\to L^\infty}\le Ch^{-n/2}\left|{\rm Im}\,z\right|^{-q},\eqno{(2.14)}$$
for $z\in{\rm supp}\,\widetilde\varphi$, ${\rm Im}\,z\neq 0$, with constants $C,q>0$ independent of $z$ and $h$.
On the other hand, using (2.6), (2.9), (2.10) and the resolvent identity, we get
 $$\left\|\left(\varphi(P(h))-\varphi(P_0(h))\right)\langle x\rangle^{s}\right\|_{L^2\to L^\infty}$$ 
$$\le Ch^{\nu}\sum_{|\alpha|\le 2}
\int_{\bf C}\left|\frac{\partial\widetilde\varphi}{\partial\bar z}(z)\right|\left\|\left(P_0(h)-z\right)^{-1}
\right\|_{L^2\to L^\infty}\left\|\langle x\rangle^{-\delta}{\cal D}_x^\alpha\left(P(h)-z\right)^{-1}
\langle x\rangle^{s}\right\|_{L^2\to L^2}L(dz).\eqno{(2.15)}$$
In this case (2.5) follows from (2.7), (2.12), (2.14) and (2.15). Let now $n\ge 4$. Then (2.14) is no longer 
true because the singularity of the kernel of the free resolvent on the diagonal gets too strong. 
In this case we will derive (2.5) from the following

\begin{lemma} Given any $0<\varepsilon\ll 1$, the free resolvent can be decomposed as
$$\left(P_0(h)-z\right)^{-1}=\sum_{j=1}^3{\cal B}_\varepsilon^{(j)}(z,h),$$
where ${\cal B}_\varepsilon^{(j)}(z,h)$, $j=1,3$, are analytic on ${\rm supp}\,\widetilde\varphi$. 
Moreover, for $s\ge 0$, $z\in{\rm supp}\,\widetilde\varphi$, we have the estimates
$$\left\|{\cal D}_x^\alpha{\cal B}_\varepsilon^{(1)}(z,h)\right\|_{L^1\to L^1}\le 
C\varepsilon^{1-\frac{|\alpha|}{2}},\quad |\alpha|\le 2,\eqno{(2.16)}$$
$$\left\|\langle x\rangle^s{\cal B}_\varepsilon^{(2)}(z,h)\langle x\rangle^{-s}\right\|_{L^1\to L^2}
\le C_\varepsilon h^{-\frac{n}{2}}\left|{\rm Im}\,z\right|^{-q},\eqno{(2.17)}$$
$$\left\|{\cal B}_\varepsilon^{(3)}(z,h)\right\|_{L^1\to L^2}\le C_\varepsilon h^{-\frac{n}{2}},\eqno{(2.18)}$$
with constants $C,q>0$ independent of $z$, $h$ and $\varepsilon$, and a constant $C_\varepsilon>0$ independent of $z$ and $h$.
\end{lemma}

Note that by assumption we have $r_\alpha=O(\gamma)$ for 
$|\alpha|=2$. 
It follows from (2.16) together with (2.9) and (2.10) that
$$\left\|\left(P(h)-P_0(h)\right){\cal B}_\varepsilon^{(1)}(z,h)\right\|_{L^1\to L^1}\le 
C(\varepsilon^{\frac{1}{2}}+\gamma),\eqno{(2.19)}$$
$$\left\|\langle x\rangle^\delta\left(P(h)-P_0(h)\right){\cal B}_\varepsilon^{(1)}(z,h)
\right\|_{L^1\to L^1}\le Ch^\nu,\eqno{(2.20)}$$
with a constant $C>0$ independent of $z$, $h$ and $\varepsilon$. Clearly, (2.19) implies that the operator 
$1+\left(P(h)-P_0(h)\right){\cal B}_\varepsilon^{(1)}(z,h)$ is invertible on $L^1$, provided 
$\varepsilon,\gamma>0$ are taken small enough, independent of $h$, with an inverse analytic 
on ${\rm supp}\,\widetilde\varphi$. Therefore, we can write
$$\left(P(h)-z\right)^{-1}-\left(P_0(h)-z\right)^{-1}=\sum_{j=1}^4{\cal F}_j(z,h),\eqno{(2.21)}$$
where
$${\cal F}_j(z,h)=-{\cal B}_\varepsilon^{(j)}(z,h)\left(P(h)-P_0(h)\right){\cal B}_\varepsilon^{(1)}(z,h)
\left(1+\left(P(h)-P_0(h)\right){\cal B}_\varepsilon^{(1)}(z,h)\right)^{-1},$$
$j=1,2,3,$ and
$${\cal F}_4(z,h)=$$ $$-\left(P(h)-z\right)^{-1}\left(P(h)-P_0(h)\right)\left({\cal B}_\varepsilon^{(2)}(z,h)
+{\cal B}_\varepsilon^{(3)}(z,h)\right)\left(1+\left(P(h)-P_0(h)\right){\cal B}_\varepsilon^{(1)}(z,h)\right)^{-1}.$$
Clearly, ${\cal F}_j(z,h)$, $j=1,3$, are analytic on ${\rm supp}\,\widetilde\varphi$, so in view of (2.21) we can write
$$\varphi(P(h))-\varphi(P_0(h))=\sum_{j=2,4}\frac{1}{\pi}\int_{\bf C}\frac{\partial\widetilde\varphi}{\partial\bar z}(z)
{\cal F}_j(z,h)L(dz).\eqno{(2.22)}$$
By (2.17) and (2.20),
$$\left\|\langle x\rangle^s{\cal F}_2(z,h)\right\|_{L^1\to L^2}\le Ch^\nu \left\|\langle x\rangle^s{\cal B}_\varepsilon^{(2)}(z,h)
\langle x\rangle^{-\delta}\right\|_{L^1\to L^2}\le C_\varepsilon h^{\nu-\frac{n}{2}}\left|{\rm Im}\,z\right|^{-q_1}.\eqno{(2.23)}$$
By (2.9), (2.10), (2.12), (2.17) and (2.18),
$$\left\|\langle x\rangle^s{\cal F}_4(z,h)\right\|_{L^1\to L^2}$$ 
$$\le C\left\|\langle x\rangle^s(P(h)-z)^{-1}(P(h)-P_0(h))\right\|_{L^2\to L^2}
\left(\left\|{\cal B}_\varepsilon^{(2)}(z,h)\right\|_{L^1\to L^2}+\left\|{\cal B}_\varepsilon^{(3)}(z,h)
\right\|_{L^1\to L^2}\right)$$
$$\le C_\varepsilon h^{\nu-\frac{n}{2}}\left|{\rm Im}\,z\right|^{-q_2}.\eqno{(2.24)}$$
By (2.7), (2.22), (2.23) and (2.24), we conclude
$$\left\|\langle x\rangle^s(\varphi(P(h))-\varphi(P_0(h)))\right\|_{L^1\to L^2}\le C_\varepsilon h^{\nu-\frac{n}{2}},$$
which is equivalent to (2.5).
\eproof

{\it Proof of Lemma 2.2.} Let $\phi\in C_0^\infty([1,2])$ be such that $\int\phi(\theta)d\theta=1$. 
Given any $0<\varepsilon\ll 1$, write
$[0,+\infty)=\cup_{j=1}^3I_j(\varepsilon)$, where $I_1(\varepsilon)=[0,\varepsilon]$, 
$I_2(\varepsilon)=[\varepsilon,\varepsilon^{-1}]$, $I_3(\varepsilon)=[\varepsilon^{-1},+\infty)$. Set
$$\chi_\varepsilon^{(j)}(\sigma)=\sigma\int_{I_j(\varepsilon)}\phi(\sigma\theta)d\theta,$$
$${\cal B}_\varepsilon^{(j)}(z,h)=\left(P_0(h)-z\right)^{-1}\chi_\varepsilon^{(j)}(P_0(h))=
\int_{I_j(\varepsilon)}\psi(\theta P_0(h),\theta z)d\theta,$$
where
$$\psi(\lambda,w)=\lambda(\lambda-w)^{-1}\phi(\lambda).$$
Since ${\rm supp}\,\widetilde\varphi$ is a compact disjoint from zero, taking $\varepsilon>0$ small enough, 
we can arrange that $\theta z$ does not belong to the support of $\phi$ as long as 
$\theta\in I_1(\varepsilon)\cup I_3(\varepsilon)$ and $z\in {\rm supp}\,\widetilde\varphi$. 
Therefore, the operator-valued functions ${\cal B}_\varepsilon^{(j)}(\cdot,h)$, $j=1,3$, 
are analytic on ${\rm supp}\,\widetilde\varphi$.
We also have the bounds
$$\left|\partial_\lambda^k\left(\lambda^{|\alpha|/2}\psi(\lambda,\theta z)\right)\right|
\le C_k,\quad \theta\in I_1(\varepsilon),\eqno{(2.25)}$$
$$\left|\partial_\lambda^k\psi(\lambda,\theta z)\right|\le C_k\left|{\rm Im}\,z\right|^{-k-1},
\quad \theta\in I_2(\varepsilon),\eqno{(2.26)}$$
$$\left|\partial_\lambda^k\psi(\lambda,\theta z)\right|\le C_k\theta^{-1},\quad \theta\in I_3(\varepsilon),\eqno{(2.27)}$$
for $z\in{\rm supp}\,\widetilde\varphi$ and all integers $k\ge 0$. Recall now that given any function 
$f\in C_0^\infty({\bf R})$ and any $h>0$,
the operator $f(P_0(h))$ satisfies the estimates (e.g. see Lemma A.1 of \cite{kn:MV})
$$\left\|f(P_0(h))\right\|_{L^1\to L^1}\le \widetilde C\sum_{k=0}^N\sup\left|\partial_\lambda^k f(\lambda)\right|,\eqno{(2.28)}$$
$$\left\|\langle x\rangle^sf(P_0(h))\langle x\rangle^{-s}\right\|_{L^1\to L^2}\le \widetilde Ch^{-n/2}
\langle h\rangle^{|s|}\sum_{k=0}^{N_s}\sup\left|\partial_\lambda^k f(\lambda)\right|,\quad\forall s\in{\bf R},\eqno{(2.29)}$$
where $N$ and $N_s$ are integers independent of $f$ and $h$, while $\widetilde C>0$ 
is a constant depending only on the support of $f$. If $|\alpha|\le 1$, by (2.25) and (2.28), we get
$$\left\|{\cal D}_x^\alpha{\cal B}_\varepsilon^{(1)}(z,h)\right\|_{L^1\to L^1}\le 
C\left\|P_0(h)^{|\alpha|/2}{\cal B}_\varepsilon^{(1)}(z,h)\right\|_{L^1\to L^1}$$ 
$$\le C\int_0^\varepsilon\left\|P_0(h)^{|\alpha|/2}\psi(\theta P_0(h),\theta z)
\right\|_{L^1\to L^1}d\theta \le C\int_0^\varepsilon \theta^{-|\alpha|/2}d\theta\le C\varepsilon^{1-|\alpha|/2}.$$
Using (2.29) together with (2.26) and (2.27), we also get 
$$\left\|\langle x\rangle^s{\cal B}_\varepsilon^{(2)}(z,h)\langle x\rangle^{-s}\right\|_{L^1\to L^2}
\le \int_\varepsilon^{\varepsilon^{-1}}\left\|\langle x\rangle^s\psi(\theta P_0(h),\theta z)
\langle x\rangle^{-s}\right\|_{L^1\to L^2}d\theta$$
$$\le Ch^{-n/2}\left|{\rm Im}\,z\right|^{-N_s-1}\int_\varepsilon^{\varepsilon^{-1}}
\theta^{-n/4}(1+\theta)^sd\theta\le C_\varepsilon h^{-n/2}\left|{\rm Im}\,z\right|^{-N_s-1},$$
 $$\left\|{\cal B}_\varepsilon^{(3)}(z,h)\right\|_{L^1\to L^2}\le \int_{\varepsilon^{-1}}^\infty
\left\|\psi(\theta P_0(h),\theta z)\right\|_{L^1\to L^2}d\theta$$
$$\le Ch^{-n/2}\int_{\varepsilon^{-1}}^\infty\theta^{-1-n/4}d\theta\le Ch^{-n/2}.$$
It remains to prove (2.16) for $|\alpha|=2$. Clearly, it suffices to show that the operator 
$\chi_\varepsilon^{(1)}(P_0(h))$ is bounded on $L^1$
uniformly in $\varepsilon$ and $h$. Since $\chi_\varepsilon^{(1)}(\sigma)=\chi_1^{(1)}(\varepsilon\sigma)$, 
we need to show that the operator $\chi_1^{(1)}(-\varepsilon h^2\Delta)$ is bounded on $L^1$ uniformly in 
$\varepsilon$ and $h$. To see this observe that the kernel of $\chi_1^{(1)}(-\varepsilon h^2\Delta)$ is of the form 
$(\varepsilon^{1/2}h)^{-n}K(|x-y|/\varepsilon^{1/2}h)$, where $K(|x-y|)$ is the kernel of 
$\chi_1^{(1)}(-\Delta)$. Hence $\chi_1^{(1)}(-\varepsilon h^2\Delta)$ is bounded on $L^1$ if and only if so is 
$\chi_1^{(1)}(-\Delta)$ and the norms
coincide. On the other hand, we have $\chi_1^{(1)}\in C^\infty({\bf R})$, $\chi_1^{(1)}(\sigma)=0$ for $\sigma\le 1$, 
$\chi_1^{(1)}(\sigma)=1$ for $\sigma\ge 2$, which implies that $\chi_1^{(1)}(-\Delta)$ is bounded on $L^1$.
\eproof

\section{Uniform estimates on weighted $L^2$ spaces}

We will prove the following

\begin{Theorem} Assume (1.3), (1.4) and (1.6) fulfilled. Let $0\le s<\delta-1$, $0<\epsilon\ll 1$. Then, we have the estimates
$$\int_{-\infty}^\infty\langle t\rangle^{2s}\left\|\langle x\rangle^{-1/2-s-\epsilon}e^{it\sqrt{P(h)}}\varphi(P(h))
\langle x\rangle^{-1/2-s-\epsilon}f\right\|_{L^2}^2dt\le C_\varepsilon\mu(h)^{2+2s+2\varepsilon}
\left\|f\right\|_{L^2}^2,\eqno{(3.1)}$$
$$\left\|\langle x\rangle^{-1/2-s-\epsilon}e^{it\sqrt{P(h)}}\varphi(P(h))\langle x\rangle^{-1/2-s-\epsilon}
\right\|_{L^2\to L^2}\le C_\varepsilon\mu(h)^{1+s+\varepsilon}\langle t\rangle^{-s},\quad \forall t,\eqno{(3.2)}$$
$$\int_{-\infty}^\infty\langle t\rangle^{2s}\left\|\langle x\rangle^{-1/2-s-\epsilon}e^{itP(h)}\varphi(P(h))
\langle x\rangle^{-1/2-s-\epsilon}f\right\|_{L^2}^2dt\le C_\varepsilon\mu(h)^{2+2s+2\varepsilon}
\left\|f\right\|_{L^2}^2,\eqno{(3.3)}$$
$$\left\|\langle x\rangle^{-1/2-s-\epsilon}e^{itP(h)}\varphi(P(h))\langle x\rangle^{-1/2-s-\epsilon}
\right\|_{L^2\to L^2}\le C_\varepsilon\mu(h)^{1+s+\varepsilon}\langle t\rangle^{-s},\quad \forall t,\eqno{(3.4)}$$
for every $0<\varepsilon\ll 1$. 
\end{Theorem}

{\it Proof.} Let us first see that (3.2) follows from (3.1). Given any $f\in L^2$, set
$$u(x,t)=\langle x\rangle^{-1/2-s-\epsilon}e^{it\sqrt{P(h)}}\varphi(P(h))\langle x\rangle^{-1/2-s-\epsilon}f.$$
It follows from (3.1) that there exists a sequence $t_k\to \infty$ such that
$$\lim_{t_k\to \infty}\left\|u(\cdot,t_k)\right\|_{L^2}=0.\eqno{(3.5)}$$
Let $\varphi_1\in C_0^\infty((E-\varepsilon_0,E+\varepsilon_0))$, $\varphi_1=1$ on supp$\,\varphi$. Using (2.2) we have
$$\left|\frac{d}{dt}\left\|u(\cdot,t)\right\|^2_{L^2}\right|=2\left|{\rm Re}\,
\left\langle\partial_tu(\cdot,t),u(\cdot,t)\right\rangle_{L^2}\right|$$ $$=
2\left|{\rm Im}\,\left\langle\langle x\rangle^{-1/2-s-\epsilon}\sqrt{P(h)}\varphi_1(P(h))
\langle x\rangle^{1/2+s+\epsilon}u(x,t),u(x,t)\right\rangle_{L^2}\right|\le C\left\|u(\cdot,t)\right\|^2_{L^2}$$
with a constant $C>0$ independent of $h$ and $t$. Hence given any $t>0$, we get
$$\left\|u(\cdot,t)\right\|^2_{L^2}\le \left\|u(\cdot,t_k)\right\|^2_{L^2}
+C\int_t^{t_k}\left\|u(\cdot,\tau)\right\|^2_{L^2}d\tau,$$
which together with (3.5) imply
$$\left\|u(\cdot,t)\right\|^2_{L^2}\le C\int_t^\infty\left\|u(\cdot,\tau)\right\|^2_{L^2}d\tau.\eqno{(3.6)}$$
By (3.6)
$$t^{2s}\left\|u(\cdot,t)\right\|^2_{L^2}\le C\int_0^\infty\tau^{2s}\left\|u(\cdot,\tau)\right\|^2_{L^2}
d\tau\le C_\varepsilon\mu(h)^{2s+2+2\varepsilon}\|f\|^2_{L^2},$$
which is the desired bound. 
The fact that (3.3) implies (3.4) can be proved in precisely the same way. We will next derive (3.1) and (3.3) from the following

\begin{prop} Assume (1.3), (1.4) and (1.6) fulfilled. Let $0\le s<\delta-1$, $0<\epsilon\ll 1$. Then, 
$${\cal D}_x^\alpha R^\pm_{1/2+s+\epsilon}(\cdot,h)\in C^s\left((E-\varepsilon_0,E+\varepsilon_0);{\cal L}(L^2)\right)$$ and
$$\left\|{\cal D}_x^\alpha R^\pm_{1/2+s+\epsilon}(\cdot,h)\right\|_{C^s}\le C\mu(h)^{1+s},\eqno{(3.7)}$$
where $0\le |\alpha|\le 2$.
\end{prop}

Observe first that it suffices to bound the integral in the left-hand side of (3.1) over the interval $[1,\infty)$ only, 
since over $(-\infty,1]$
it can be treated similarly while over $[-1,1]$ it is trivial. Thus, it suffices to prove the bound
$$\int_{2^k}^{2^{k+1}}\left\|u(\cdot,t)\right\|^2_{L^2}dt\le C 2^{-2k(s+\varepsilon)}\mu(h)^{2+2s+2\varepsilon}
\left\|f\right\|^2_{L^2},\eqno{(3.8)}$$
for every integer $k\ge 0$ and every $0\le\varepsilon\ll 1$. Let $\phi\in C^\infty({\bf R})$, $0\le\phi\le 1$, 
$\phi(t)=0$ for $t\le 1/3$, $\phi(t)=1$ for $t\ge 1/2$. We have
$$\left(\partial_t^2+P(h)\right)\langle x\rangle^{1/2+s+\epsilon}\phi(t)u(x,t)=2i\phi'(t)\sqrt{P(h)}
\langle x\rangle^{1/2+s+\epsilon}u(x,t)+
\phi''(t)\langle x\rangle^{1/2+s+\epsilon}u(x,t)$$ $$=:\varphi_1(P(h))\langle x\rangle^{-1/2-s-\epsilon}v(x,t).$$
In view of (2.2) we have
$$\left\|v(\cdot,t)\right\|_{L^2}\le C\|f\|_{L^2},\quad\forall t,\eqno{(3.9)}$$
with a contsant $C>0$ independent of $t$ and $h$. By Duhamel's formula we get
$$\phi(t)u(x,t)=\int_0^t\langle x\rangle^{-1/2-s-\epsilon}\frac{\sin\left((t-\tau)
\sqrt{P(h)}\right)}{\sqrt{P(h)}}\varphi_1(P(h))\langle x\rangle^{-1/2-s-\epsilon}v(x,\tau)d\tau.\eqno{(3.10)}$$
Taking the Fourier transform with respect to $t$, we deduce from (3.10)
$$\widehat{\phi u}(x,\lambda)=T(\lambda,h)\widehat{v}(x,\lambda),\quad\lambda\in{\bf R},\eqno{(3.11)}$$
where $$T(\lambda,h)=\langle x\rangle^{-1/2-s-\epsilon}\left(P(h)-\lambda^2+i0\right)^{-1}\varphi_1(P(h))
\langle x\rangle^{-1/2-s-\epsilon}.$$
It is easy to see that $T(\cdot,h)\in C^s\left({\bf R};{\cal L}(L^2)\right)$ and
$$\left\|T(\cdot,h)\right\|_{C^s}\le C\mu(h)^{1+s}.\eqno{(3.12)}$$
Indeed, if $\lambda^2$ belongs to a small neighbourhood, $K$, of supp$\,\varphi_1$, $K\subset (E-\varepsilon_0,E+\varepsilon_0)$, 
then this follows from (2.2) and Proposition 3.2. Let $\lambda^2\in {\bf R}\setminus K$. Then, for every integer $m\ge 0$, we have
$$\left\|\frac{d^m}{d\lambda^m}\left(P(h)-\lambda^2\right)^{-1}\varphi_1(P(h))\right\|_{L^2\to L^2}
\le C_m\sum_{j=0}^m\langle\lambda\rangle^j
\left\|\left(P(h)-\lambda^2\right)^{-1-j}\varphi_1(P(h))\right\|_{L^2\to L^2}$$
$$\le C_m\sum_{j=0}^m\langle\lambda\rangle^j\sup_{y\in{\rm supp}\,\varphi_1}|y-\lambda^2|^{-1-j}\le 
C'_m \sum_{j=0}^m\langle\lambda\rangle^{-2-j}\le Const.$$
Let now $\rho\in C_0^\infty([1/3,1/2])$, $\rho\ge 0$, such that $\int\rho(\sigma)d\sigma=1$. Given a parameter $0<\theta\le 1$, set
$$T_\theta(\lambda,h)=\theta^{-1}\int T(\lambda+\sigma,h)\rho(\sigma/\theta)d\sigma.$$
It follows from (3.12) that the operator-valued function $T_\theta(\lambda,h)$ satisfies the bounds
$$\left\|\partial_\lambda^jT_\theta(\lambda,h)\right\|_{L^2\to L^2}\le C\mu(h)^{1+j},\quad 0\le j\le [s],\eqno{(3.13)}$$
$$\left\|\partial_\lambda^j\left(T-T_\theta\right)(\lambda,h)\right\|_{L^2\to L^2}
\le \theta^{-1}\int\left\|\partial_\lambda^jT(\lambda+\sigma,h)-\partial_\lambda^jT(\lambda,h)
\right\|_{L^2\to L^2}\rho(\sigma/\theta)d\sigma
$$ $$\le C\mu(h)^{1+s}\theta^{s-[s]},\quad 0\le j\le [s],\eqno{(3.14)}$$
$$\left\|\partial_\lambda^jT_\theta(\lambda,h)\right\|_{L^2\to L^2}\le \theta^{-2}
\int\left\|\partial_\lambda^{j-1}T(\lambda+\sigma,h)-\partial_\lambda^{j-1}T(\lambda,h)
\right\|_{L^2\to L^2}\left|\rho'(\sigma/\theta)\right|d\sigma
$$ $$\le C\mu(h)^{1+s}\theta^{s-[s]-1},\quad j=[s]+1.\eqno{(3.15)}$$
Define the function $u_\theta(t,x)$ by the relation
$$\widehat{u_\theta}(x,\lambda)=T_\theta(\lambda,h)\widehat{v}(x,\lambda).$$
Using (3.9), (3.14) together with Plancherel's identity, we obtain
$$\int_{-\infty}^\infty|t|^{2[s]}\left\|\phi u(\cdot,t)-u_\theta(\cdot,t)\right\|^2_{L^2}dt=\int_{-\infty}^\infty
\left\|\partial_\lambda^{[s]}\left(\widehat{\phi u}(\cdot,\lambda)-\widehat{u_\theta}(\cdot,\lambda)\right)
\right\|^2_{L^2}d\lambda$$
$$\le C\sum_{j=0}^{[s]}\int_{-\infty}^\infty\left\|\partial_\lambda^{j}\left(T-T_\theta\right)
(\lambda,h)\partial_\lambda^{[s]-j} \widehat{v}(\cdot,\lambda)\right\|^2_{L^2}d\lambda$$ 
$$\le C\mu(h)^{2+2s}\theta^{2s-2[s]}\sum_{j=0}^{[s]}\int_{-\infty}^\infty\left\|\partial_\lambda^{[s]-j} 
\widehat{v}(\cdot,\lambda)\right\|^2_{L^2}d\lambda$$
$$= C\mu(h)^{2+2s}\theta^{2s-2[s]}\sum_{j=0}^{[s]}\int_{1/3}^{1/2} t^{2[s]-2j}\left\|v(\cdot,t)
\right\|^2_{L^2}dt\le C\mu(h)^{2+2s}\theta^{2s-2[s]}\|f\|^2_{L^2}.$$
Hence, given a parameter $M\ge 1$, we get
$$\int_{M}^{2M}\left\|u(\cdot,t)-u_\theta(\cdot,t)\right\|^2_{L^2}dt\le C\mu(h)^{2+2s}
\theta^{2s-2[s]}M^{-2[s]}\|f\|^2_{L^2}.\eqno{(3.16)}$$
Similarly, using (3.9), (3.13), (3.15) together with Plancherel's identity, we obtain
$$\int_{-\infty}^\infty|t|^{2[s]+2}\left\|u_\theta(\cdot,t)\right\|^2_{L^2}dt=\int_{-\infty}^\infty
\left\|\partial_\lambda^{[s]+1}\widehat{u_\theta}(\cdot,\lambda)\right\|^2_{L^2}d\lambda$$
$$\le C\sum_{j=0}^{[s]+1}\int_{-\infty}^\infty\left\|\partial_\lambda^{j}T_\theta(\lambda,h)
\partial_\lambda^{[s]+1-j} \widehat{v}(\cdot,\lambda)\right\|^2_{L^2}d\lambda$$ $$\le C\mu(h)^{2+2s}
\theta^{2s-2[s]-2}\sum_{j=0}^{[s]+1}\int_{-\infty}^\infty\left\|\partial_\lambda^{[s]+1-j} 
\widehat{v}(\cdot,\lambda)\right\|^2_{L^2}d\lambda$$
$$= C\mu(h)^{2+2s}\theta^{2s-2[s]-2}\sum_{j=0}^{[s]+1}\int_{1/3}^{1/2} t^{2[s]+2-2j}
\left\|v(\cdot,t)\right\|^2_{L^2}dt\le C\mu(h)^{2+2s}\theta^{2s-2[s]-2}\|f\|^2_{L^2},$$
which implies
$$\int_{M}^{2M}\left\|u_\theta(\cdot,t)\right\|^2_{L^2}dt\le C\mu(h)^{2+2s}\theta^{2s-2[s]-2}M^{-2[s]-2}\|f\|^2_{L^2}.\eqno{(3.17)}$$
Taking $\theta=M^{-1}$ we deduce from (3.16) and (3.17)
$$\int_{M}^{2M}\left\|u(\cdot,t)\right\|^2_{L^2}dt\le C\mu(h)^{2+2s}M^{-2s}\|f\|^2_{L^2}.\eqno{(3.18)}$$
Observe finally that the estimates (3.13)-(3.15) hold true with $s$ replaced by $s+\varepsilon$, $\forall 0\le\varepsilon\ll 1$, 
and hence so does (3.18), which in turn proves (3.8).

The estimate (3.3) can be proved in the same way. The only difference is that the function
$$w(x,t)=\langle x\rangle^{-1/2-s-\epsilon}e^{itP(h)}\varphi(P(h))\langle x\rangle^{-1/2-s-\epsilon}f$$
satisfies the identity
$$\widehat{\phi w}(x,\lambda)=\widetilde T(\lambda,h)\widehat{v}(x,\lambda),\quad\lambda\in{\bf R},$$
where $$\widetilde T(\lambda,h)=\langle x\rangle^{-1/2-s-\epsilon}\left(P(h)-\lambda+i0\right)^{-1}\varphi_1(P(h))
\langle x\rangle^{-1/2-s-\epsilon}$$
belongs again to $C^s\left({\bf R};{\cal L}(L^2)\right)$ and satisfies (3.12), while the function $v(x,t)$ 
is compactly supported in $t$ and satisfies (3.9).
\eproof

{\it Proof of Proposition 3.2.} We will use the commutator identity
$$\Delta+\frac{1}{2}\left[x\cdot\nabla,\Delta\right]=0,$$
which we rewrite as follows
$$P(h)+\frac{1}{2}\left[x\cdot\nabla,P(h)\right]=P(h)-P_0(h)+\frac{1}{2}\left[x\cdot\nabla,P(h)-P_0(h)\right]=:{\cal Q}(h).
\eqno{(3.19)}$$
Given any $z\in{\bf C}$, Im$\,z\neq 0$, we deduce from (3.19)
$$P(h)-z+\frac{1}{2}\left[x\cdot\nabla,P(h)-z\right]=-z+{\cal Q}(h),$$
which yields the identity
$$\left(P(h)-z\right)^{-1}-\frac{1}{2}\left[x\cdot\nabla,\left(P(h)-z\right)^{-1}\right]$$ $$=-z\left(P(h)-z\right)^{-2}+
\left(P(h)-z\right)^{-1}{\cal Q}(h)\left(P(h)-z\right)^{-1}.\eqno{(3.20)}$$
We will first consider the case $s=m$, where $0\le m<\delta-1$ is an integer. We will proceed by induction. 
When $m=0$ the assertion is true by assumption. Suppose it is true for all integers $0\le k\le m-1$. 
We differentiate $m-1$ times the identity (3.20) with respect to $z$ to get
$$z\left(P(h)-z\right)^{-m-1}=\widetilde c_m\left(P(h)-z\right)^{-m}+\frac{1}{2}\left[x\cdot\nabla,\left(P(h)-z\right)^{-m}\right]$$
$$+\sum_{k=1}^mc_k\left(P(h)-z\right)^{-k}{\cal Q}(h)\left(P(h)-z\right)^{k-m-1},\eqno{(3.21)}$$
which in turn leads to the identity
$$z\frac{d^{m}}{dz^m}R_{1/2+m+\epsilon}(z,h)=\widetilde c_m\langle x\rangle^{-1}\frac{d^{m-1}}{dz^{m-1}}
R_{-1/2+m+\epsilon}(z,h)\langle x\rangle^{-1}$$ $$+\frac{1}{2}\langle x\rangle^{-1/2-m-\epsilon}x\cdot\nabla
\langle x\rangle^{-1/2+m+\epsilon}\frac{d^{m-1}}{dz^{m-1}}R_{-1/2+m+\epsilon}(z,h)\langle x\rangle^{-1}$$ 
$$-\frac{1}{2}\langle x\rangle^{-1}
\frac{d^{m-1}}{dz^{m-1}}R_{-1/2+m+\epsilon}(z,h)\langle x\rangle^{-1/2+m+\epsilon}x\cdot\nabla\langle x\rangle^{-1/2-m-\epsilon}$$
$$+\sum_{k=1}^mc_k\langle x\rangle^{k-m-1}\frac{d^{k-1}}{dz^{k-1}}R_{-1/2+k+\epsilon}(z,h)\widetilde{\cal Q}_k(h)
\frac{d^{m-k}}{dz^{m-k}}R_{1/2+m-k+\epsilon}(z,h)\langle x\rangle^{-k},\eqno{(3.22)}$$
where
$$\widetilde{\cal Q}_k(h)=\langle x\rangle^{-1/2+k+\epsilon}{\cal Q}(h)\langle x\rangle^{1/2+m-k+\epsilon}.$$
A simple computation shows that
$${\cal Q}(h)=\frac{1}{2}\sum_{j=1}^n\left(b_j(x,h){\cal D}_{x_j}+{\cal D}_{x_j}b_j(x,h)\right)+V(x,h)+\frac{1}{2}
\left[x\cdot h\nabla,V^1(x,h)\right]$$ $$+\frac{1}{2}\sum_{i,j=1}^n{\cal D}_{x_i}x\cdot \nabla a_{ij}(x,h){\cal D}_{x_j}+
\frac{1}{2}\sum_{j=1}^n\left(x\cdot\nabla b_j(x,h){\cal D}_{x_j}+{\cal D}_{x_j}x\cdot\nabla b_j(x,h)\right)+\frac{1}{2}x\cdot
\nabla V^0(x).$$
Hence, in view of (1.4), we have that the operators $\widetilde{\cal Q}_k(h)$ are of the form
$$\widetilde{\cal Q}_k(h)=\sum_{\alpha,\beta\in\Omega}{\cal D}_{x}^\alpha q_{\alpha,\beta}^{(k)}(x,h){\cal D}_{x}^\beta,$$
where $\Omega$ is the set of all multi-indices such that $0\le|\alpha|\le 2$, $0\le|\beta|\le 2$, $|\alpha|+|\beta|\le 3$, 
and the coefficients satisfy
$$\left|q_{\alpha,\beta}^{(k)}(x,h)\right|\le C,\eqno{(3.23)}$$
with a constant $C>0$ independent of $x$ and $h$. By (3.22) and (3.23) we obtain
$$|z|\left\|\frac{d^{m}}{dz^m}R_{1/2+m+\epsilon}(z,h)\right\|_{L^2\to L^2}\le Ch^{-1}\sum_{0\le|\alpha|\le 1}
\left\|{\cal D}_{x}^\alpha\frac{d^{m-1}}{dz^{m-1}}R_{-1/2+m+\epsilon}(z,h)\right\|_{L^2\to L^2}$$ 
$$+Ch^{-1}\sum_{0\le|\alpha|\le 1}\left\|{\cal D}_{x}^\alpha\frac{d^{m-1}}{dz^{m-1}}R_{-1/2+m+\epsilon}(\bar z,h)
\right\|_{L^2\to L^2}$$
$$+C\sum_{k=1}^m\sum_{\alpha,\beta\in\Omega}\left\|{\cal D}_{x}^\alpha\frac{d^{k-1}}{dz^{k-1}}R_{-1/2+k+\epsilon}
(\bar z,h)\right\|_{L^2\to L^2}\left\|{\cal D}_{x}^\beta\frac{d^{m-k}}{dz^{m-k}}R_{1/2+m-k+\epsilon}(z,h)
\right\|_{L^2\to L^2},\eqno{(3.24)}$$
with a constant $C>0$ independent of $z$ and $h$. Applying (3.24) with $z$ replaced by $z\pm i\varepsilon$, 
$z\in[E-\varepsilon_0,E+\varepsilon_0]$, $0<\varepsilon\ll 1$, and taking the limit as $\varepsilon\to 0$, we get
$$\left\|\frac{d^{m}}{dz^m}R_{1/2+m+\epsilon}^\pm(z,h)\right\|_{L^2\to L^2}\le Ch^{-1}\sum_\pm\sum_{0\le|\alpha|\le 1}
\left\|{\cal D}_{x}^\alpha\frac{d^{m-1}}{dz^{m-1}}R_{-1/2+m+\epsilon}^\pm(z,h)\right\|_{L^2\to L^2}$$ 
$$+C\sum_{k=1}^m\sum_{\alpha,\beta\in\Omega}\left\|{\cal D}_{x}^\alpha\frac{d^{k-1}}{dz^{k-1}}R_{-1/2+k+\epsilon}^\mp(z,h)
\right\|_{L^2\to L^2}\left\|{\cal D}_{x}^\beta\frac{d^{m-k}}{dz^{m-k}}R_{1/2+m-k+\epsilon}^\pm(z,h)\right\|_{L^2\to L^2}
\le C\mu(h)^{m+1},$$
provided (3.7) holds for all integers $s\le m-1$. Thus we get (3.7) with $s=m$ and $\alpha=0$. 
The fact that it holds for all multi-indices
$|\alpha|\le 2$ follows from the ellipticity condition (1.3). 

Let now $s=m+\nu$, where $0<\nu<1$ and $m$ is an integer such that $0\le m<\delta-1-\nu$. In this case it suffices to show that
$$\left\|\frac{d^{m+1}}{dz^{m+1}}R_{1/2+m+\nu+\epsilon}(z\pm i\varepsilon,h)\right\|_{L^2\to L^2}\le C\mu(h)^{m+\nu+1}
\varepsilon^{-1+\nu},\quad\forall z\in[E-\varepsilon_0,E+\varepsilon_0]. \eqno{(3.25)}$$
Indeed, (3.25) implies 
$$\left\|\frac{d^{m}}{dz^{m}}R_{1/2+m+\nu+\epsilon}(z\pm i\varepsilon,h)-\frac{d^{m}}{dz^{m}}R_{1/2+m+\nu+\epsilon}^\pm(z,h)
\right\|_{L^2\to L^2}$$
$$\le \int_0^\varepsilon\left\|\frac{d^{m+1}}{dz^{m+1}}R_{1/2+m+\nu+\epsilon}(z\pm i\sigma,h)\right\|_{L^2\to L^2}d\sigma$$
$$\le C\mu(h)^{m+\nu+1}\int_0^\varepsilon\sigma^{-1+\nu}d\sigma\le 
C\mu(h)^{m+\nu+1}\varepsilon^{\nu}. \eqno{(3.26)}$$
Now, given any $z_1,z_2\in[E-\varepsilon_0,E+\varepsilon_0]$, $0<|z_1-z_2|\le 1$, by (3.25) and (3.26), we get
$$\left\|\frac{d^{m}}{dz^{m}}R_{1/2+m+\nu+\epsilon}^\pm(z_1,h)-\frac{d^{m}}{dz^{m}}R_{1/2+m+\nu+\epsilon}^\pm(z_2,h)
\right\|_{L^2\to L^2}$$
$$\le\sum_{j=1}^2\left\|\frac{d^{m}}{dz^{m}}R_{1/2+m+\nu+\epsilon}(z_j\pm i\varepsilon,h)-\frac{d^{m}}{dz^{m}}
R_{1/2+m+\nu+\epsilon}^\pm(z_j,h)\right\|_{L^2\to L^2}$$
$$+\left\|\frac{d^{m}}{dz^{m}}R_{1/2+m+\nu+\epsilon}(z_1\pm i\varepsilon,h)-\frac{d^{m}}{dz^{m}}
R_{1/2+m+\nu+\epsilon}(z_2\pm i\varepsilon,h)\right\|_{L^2\to L^2}$$
$$\le C\mu(h)^{m+\nu+1}\left(\varepsilon^{\nu}+|z_1-z_2|\varepsilon^{-1+\nu}\right)
\le C\mu(h)^{m+\nu+1}|z_1-z_2|^{\nu}\eqno{(3.27)}$$
if we take $\varepsilon=|z_1-z_2|$. So, in this case (3.7) with $\alpha=0$ follows from (3.27). 
For any multi-index $|\alpha|\le 2$, it follows from (1.3).

Using (3.21) and proceeding by induction as above, it is easy to see that (3.25) follows from the following

\begin{lemma}Let $z\in[E-\varepsilon_0,E+\varepsilon_0]$, $0<\epsilon,\varepsilon\ll 1$, $0\le\nu\le 1$. Then
$$\left\|\langle x\rangle^{-\nu/2-\epsilon}\left(P(h)-z\pm i\varepsilon\right)^{-1}
\langle x\rangle^{-1/2-\epsilon}\right\|_{L^2\to L^2}\le C\mu(h)^{\frac{\nu+1}{2}}\varepsilon^{\frac{\nu-1}{2}}. \eqno{(3.28)}$$
\end{lemma}

{\it Proof.} When $\nu=1$ (3.28) follows from (1.6). To prove (3.28) for $\nu=0$ we will use the identity
$$\left(P(h)-z\mp i\varepsilon\right)^{-1}\left(P(h)-z\pm i\varepsilon\right)^{-1}=\frac{\pm 2}{\varepsilon}
\left(\left(P(h)-z\pm i\varepsilon\right)^{-1}-\left(P(h)-z\mp i\varepsilon\right)^{-1}\right).$$
Hence, the operator
$${\cal A}=\left(P(h)-z\pm i\varepsilon\right)^{-1}\langle x\rangle^{-1/2-\epsilon}$$
satisfies
$$\left\|{\cal A}^*{\cal A}\right\|_{L^2\to L^2}\le 2\varepsilon^{-1}\sum_\pm
\left\|\langle x\rangle^{-1/2-\epsilon}\left(P(h)-z\pm i\varepsilon\right)^{-1}\langle x\rangle^{-1/2-\epsilon}
\right\|_{L^2\to L^2}\le C\mu(h)\varepsilon^{-1},$$
where we have also used (1.6). Let now $0<\nu<1$. Given a set ${\cal M}\subset{\bf R}^n$, 
denote by $\eta({\cal M})$ the characteristic function of
${\cal M}$. Let $M>1$ be a parameter to be fixed later on. We have
$$\left\|\langle x\rangle^{-\nu/2-\epsilon}\left(P(h)-z\pm i\varepsilon\right)^{-1}\langle x\rangle^{-1/2-\epsilon}
\right\|_{L^2\to L^2}$$ $$\le
\left\|\langle x\rangle^{-\nu/2-\epsilon}\eta(\langle x\rangle\ge M)\left(P(h)-z\pm i\varepsilon\right)^{-1}
\langle x\rangle^{-1/2-\epsilon}\right\|_{L^2\to L^2}$$
$$+
\left\|\langle x\rangle^{-\nu/2-\epsilon}\eta(\langle x\rangle\le M)\left(P(h)-z\pm i\varepsilon\right)^{-1}
\langle x\rangle^{-1/2-\epsilon}\right\|_{L^2\to L^2}$$
$$\le M^{-\nu/2}
\left\|\left(P(h)-z\pm i\varepsilon\right)^{-1}\langle x\rangle^{-1/2-\epsilon}\right\|_{L^2\to L^2}$$
$$+M^{(1-\nu)/2}
\left\|\langle x\rangle^{-1/2-\epsilon}\left(P(h)-z\pm i\varepsilon\right)^{-1}\langle x\rangle^{-1/2-\epsilon}
\right\|_{L^2\to L^2}$$
$$\le CM^{-\nu/2}\mu(h)^{1/2}\varepsilon^{-1/2}+CM^{(1-\nu)/2}\mu(h)\le C\mu(h)^{\frac{\nu+1}{2}}\varepsilon^{\frac{\nu-1}{2}},$$
if we choose $M=\mu(h)\varepsilon^{-1}$.
\eproof

\section{Dispersive estimates}

We will first prove the following

\begin{prop} Under the assumptions of Theorem 1.1, for all $t\neq 0$, $0<\epsilon, \varepsilon\ll 1$, $0\le s\le\frac{n-1}{2}$,  
we have the estimates
$$\left\|\langle x\rangle^{-\sigma}\left(e^{it\sqrt{P(h)}}\varphi(P(h))-e^{it\sqrt{P_0(h)}}\varphi(P_0(h))\right)\langle x\rangle^{-1/2-s-\sigma-\epsilon}
\right\|_{L^2\to L^\infty}$$ $$\le C_\varepsilon h^{\nu-\frac{n+1}{2}}\mu(h)^{1+s+\sigma+\varepsilon}|t|^{-s-\sigma},\eqno{(4.1)}$$
$$\left\|\langle x\rangle^{-\sigma}\left(e^{itP(h)}\varphi(P(h))-e^{itP_0(h)}\varphi(P_0(h))\right)\langle x\rangle^{-1-s-\sigma-\epsilon}
\right\|_{L^2\to L^\infty}$$ $$\le C_\varepsilon h^{\nu-\frac{n+1}{2}}\mu(h)^{3/2+s+\sigma+\varepsilon}|t|^{-s-\sigma-1/2}.\eqno{(4.2)}$$
\end{prop}

{\it Proof.} Recall first that the free groups satisfy the estimates (see the appendix)
$$\left\|\langle x\rangle^{-\sigma}e^{it\sqrt{P_0(h)}}\varphi(P_0(h))\langle x\rangle^{-1/2-s-\sigma-\epsilon}\right\|_{L^2\to L^\infty}
\le C h^{-s-\sigma-\frac{n+1}{2}}|t|^{-s-\sigma},\eqno{(4.3)}$$
$$\int_{-\infty}^\infty|t|^{2s+2\sigma}\left\|\langle x\rangle^{-1/2-s-\sigma-\epsilon}e^{it\sqrt{P_0(h)}}\varphi(P_0(h))\langle x\rangle^{-\sigma}f\right\|_{L^2}^2dt
\le Ch^{-n-1-2s-2\sigma}\|f\|_{L^1}^2,\eqno{(4.4)}$$
$$\left\|\langle x\rangle^{-\sigma}e^{itP_0(h)}\varphi(P_0(h))\langle x\rangle^{-1/2-s-\sigma-\epsilon}\right\|_{L^2\to L^\infty}
\le C h^{-s-\sigma-\frac{n+1}{2}}|t|^{-s-\sigma-1/2},\eqno{(4.5)}$$
$$\int_{-\infty}^\infty|t|^{2s+2\sigma}\left\|\langle x\rangle^{-1/2-s-\sigma-\epsilon}e^{itP_0(h)}\varphi(P_0(h))\langle x\rangle^{-\sigma}f
\right\|_{L^2}^2dt\le Ch^{-n-1-2s-2\sigma}\|f\|_{L^1}^2.\eqno{(4.6)}$$
Without loss of generality we may suppose that $t>0$. To prove (4.1) observe first that Duhamel's formula for the 
wave equation implies the identity (e.g. see Section 3 of \cite{kn:V1})
$$e^{it\sqrt{P(h)}}\varphi(P(h))-e^{it\sqrt{P_0(h)}}\varphi(P_0(h))=\Phi_1(t,h)+\Phi_2(t,h),\eqno{(4.7)}$$
where
$$\Phi_1(t,h)=\left(\varphi_1(P(h))-\varphi_1(P_0(h))\right)e^{it\sqrt{P(h)}}\varphi(P(h))$$ $$+\varphi_1(P_0(h))
\cos\left(t\sqrt{P_0(h)}\right)\left(\varphi(P(h))-\varphi(P_0(h))\right)$$ $$+i\varphi_1^\sharp(P_0(h))
\sin\left(t\sqrt{P_0(h)}\right)\left(\varphi^\sharp(P(h))-\varphi^\sharp(P_0(h))\right),$$
$$\Phi_2(t,h)=-\int_0^t\varphi_1^\sharp(P_0(h))\sin\left((t-\tau)\sqrt{P_0(h)}\right){\cal P}(h)
e^{i\tau\sqrt{P(h)}}\varphi(P(h))d\tau,$$
where
$${\cal P}(h)=\varphi_2(P_0(h))\left(P(h)-P_0(h)\right)\varphi_2(P(h)),$$
$\varphi_1,\varphi_2\in C_0^\infty([E-\varepsilon_0,E+\varepsilon_0])$, $\varphi_1=1$ on supp$\,\varphi$, 
$\varphi_2=1$ on supp$\,\varphi_1$,
$\varphi^\sharp(z)=z^{1/2}\varphi(z)$, $\varphi_1^\sharp(z)=z^{-1/2}\varphi_1(z)$. By (2.3), (2.5), (3.2) and (4.3), we get
$$\left\|\langle x\rangle^{-\sigma}\Phi_1(t,h)\langle x\rangle^{-1/2-s-\sigma-\epsilon}\right\|_{L^2\to L^\infty}$$ 
$$\le \left\|\left(\varphi_1(P(h))-\varphi_1(P_0(h))\right)\langle x\rangle^{n/2+\sigma+\epsilon}
\right\|_{L^2\to L^\infty}$$ $$\times\left\|\langle x\rangle^{-n/2-\sigma-\epsilon}e^{it\sqrt{P(h)}}
\varphi(P(h))\langle x\rangle^{-1/2-s-\sigma-\epsilon}\right\|_{L^2\to L^2}$$
$$+\left\|\langle x\rangle^{-\sigma}\cos\left(t\sqrt{P_0(h)}\right)\varphi_1(P_0(h))\langle x\rangle^{-1/2-s-\sigma-\epsilon}
\right\|_{L^2\to L^\infty}$$ $$\times\left\|\langle x\rangle^{1/2+s+\sigma+\epsilon}
\left(\varphi(P(h))-\varphi(P_0(h))\right)\langle x\rangle^{-1/2-s-\sigma-\epsilon}\right\|_{L^2\to L^2}$$ 
$$+\left\|\langle x\rangle^{-\sigma}\sin\left(t\sqrt{P_0(h)}\right)\varphi_1^\sharp(P_0(h))\langle x\rangle^{-1/2-s-\sigma-\epsilon}\right\|_{L^2\to L^\infty}$$ 
$$\times\left\|\langle x\rangle^{1/2+s+\sigma+\epsilon}\left(\varphi^\sharp(P(h))-\varphi^\sharp(P_0(h))\right)
\langle x\rangle^{-1/2-s-\sigma-\epsilon}\right\|_{L^2\to L^2}$$
$$\le C_\varepsilon h^{\nu-n/2}\mu(h)^{1+s+\sigma+\varepsilon}t^{-s-\sigma}+Ch^{\nu-s-\sigma-(n+1)/2}t^{-s-\sigma}.\eqno{(4.8)}$$
Furthermore, given any $f\in L^2$, $g\in L^1$, using (2.4), (3.1) and (4.4), we get
$$t^{s+\sigma}\left|\left\langle\Phi_2(t,h)\langle x\rangle^{-1/2-s-\sigma-\epsilon}f,\langle x\rangle^{-\sigma}g\right\rangle\right|\le C\int_0^{t/2}(t-\tau)^{s+\sigma}$$ $$\times\left|\left\langle {\cal P}(h)e^{i\tau\sqrt{P(h)}}\varphi(P(h))
\langle x\rangle^{-1/2-s-\sigma-\epsilon}f,\sin\left((t-\tau)\sqrt{P_0(h)}\right)\varphi_1^\sharp(P_0(h))\langle x\rangle^{-\sigma}g\right\rangle\right|d\tau$$
$$+ C\int_{t/2}^t\tau^{s+\sigma}\left|\left\langle {\cal P}(h)e^{i\tau\sqrt{P(h)}}\varphi(P(h))
\langle x\rangle^{-1/2-s-\sigma-\epsilon}f,\sin\left((t-\tau)\sqrt{P_0(h)}\right)\varphi_1^\sharp(P_0(h))\langle x\rangle^{-\sigma}g\right\rangle\right|d\tau$$
 $$\le Ch^\nu\int_0^{t/2}(t-\tau)^{s+\sigma}\left\|\langle x\rangle^{-1/2-\epsilon}e^{i\tau\sqrt{P(h)}}
\varphi(P(h))\langle x\rangle^{-1/2-s-\sigma-\epsilon}f\right\|_{L^2}$$ 
$$\times\left\|\langle x\rangle^{-1/2-s-\sigma-\epsilon}\sin\left((t-\tau)
\sqrt{P_0(h)}\right)\varphi_1^\sharp(P_0(h))\langle x\rangle^{-\sigma}g\right\|_{L^2}d\tau$$
$$+ Ch^\nu\int_{t/2}^t\tau^{s+\sigma}\left\|\langle x\rangle^{-1/2-s-\sigma-\epsilon}
e^{i\tau\sqrt{P(h)}}\varphi(P(h))\langle  x\rangle^{-1/2-s-\sigma-\epsilon}f\right\|_{L^2}$$ 
$$\times\left\|\langle x\rangle^{-1/2-\epsilon}\sin\left((t-\tau)\sqrt{P_0(h)}\right)\varphi_1^\sharp(P_0(h))\langle x\rangle^{-\sigma}g\right\|_{L^2}d\tau$$
 $$\le Ch^\nu\left(\int_0^\infty\left\|\langle x\rangle^{-1/2-\epsilon}e^{i\tau\sqrt{P(h)}}
\varphi(P(h))\langle x\rangle^{-1/2-s-\sigma-\epsilon}f\right\|_{L^2}^2d\tau\right)^{1/2}$$ 
$$\times\left(\int_0^\infty\tau^{2s+2\sigma}\left\|\langle x\rangle^{-1/2-s-\sigma-\epsilon}
\sin\left(\tau\sqrt{P_0(h)}\right)\varphi_1^\sharp(P_0(h))\langle x\rangle^{-\sigma}g\right\|_{L^2}^2d\tau\right)^{1/2}$$
 $$+ Ch^\nu\left(\int_0^\infty\tau^{2s+2\sigma}\left\|\langle x\rangle^{-1/2-s-\sigma-\epsilon}
e^{i\tau\sqrt{P(h)}}\varphi(P(h))\langle  x\rangle^{-1/2-s-\sigma-\epsilon}f\right\|_{L^2}^2d\tau\right)^{1/2}$$ 
$$\times\left(\int_0^\infty\left\|\langle x\rangle^{-1/2-\epsilon}
\sin\left(\tau\sqrt{P_0(h)}\right)\varphi_1^\sharp(P_0(h))\langle x\rangle^{-\sigma}g\right\|_{L^2}^2d\tau\right)^{1/2}$$
 $$\le C_\varepsilon h^{\nu-s-\sigma-\frac{n+1}{2}}\mu(h)^{1+\varepsilon}\|f\|_{L^2}\|g\|_{L^1}+ 
C_\varepsilon h^{\nu-\frac{n+1}{2}}\mu(h)^{1+s+\sigma+\varepsilon}\|f\|_{L^2}\|g\|_{L^1}.\eqno{(4.9)}$$
Clearly, (4.1) follows from (4.7), (4.8) and (4.9). To prove (4.2) we rewrite Duhamel's formula 
for the Schr\"odinger equation as follows
$$e^{itP(h)}\varphi(P(h))-e^{itP_0(h)}\varphi(P_0(h))=\Psi_1(t,h)+\Psi_2(t,h),\eqno{(4.10)}$$
where
$$\Psi_1(t,h)=\left(\varphi_1(P(h))-\varphi_1(P_0(h))\right)e^{itP(h)}\varphi(P(h))+
\varphi_1(P_0(h))e^{itP_0(h)}\left(\varphi(P(h))-\varphi(P_0(h))\right),$$
$$\Psi_2(t,h)=i\int_0^te^{i(t-\tau)P_0(h)}\varphi_1(P_0(h)){\cal P}(h)e^{i\tau P(h)}\varphi(P(h))d\tau.$$
Using (4.10) together with (4.5) and (4.6), it is easy to see that (4.2) can be proved in the same way as (4.1) above.
\eproof

Clearly, (1.8) (resp. (1.10)) follows from (4.1) and (4.3) (resp. (4.2) and (4.5)) applied with $s=\frac{n-1}{2}$.
To prove (1.9) we will use once again the identity (4.7).  By (2.5), (4.1) and (4.3), we get
$$\left\|\langle x\rangle^{-\sigma}\Phi_1(t,h)\langle x\rangle^{-\sigma}\right\|_{L^1\to L^\infty}$$ $$\le \left\|\left(\varphi_1(P(h))-\varphi_1(P_0(h))\right)
\langle x\rangle^{n/2+\sigma+\epsilon}\right\|_{L^2\to L^\infty}$$ $$\times\left\|\langle x\rangle^{-n/2-\sigma-\epsilon}
e^{it\sqrt{P(h)}}\varphi(P(h))\langle x\rangle^{-\sigma}\right\|_{L^1\to L^2}$$
$$+\left\|\langle x\rangle^{-\sigma}\cos\left(t\sqrt{P_0(h)}\right)\varphi_1(P_0(h))\langle x\rangle^{-n/2-\sigma-\epsilon}
\right\|_{L^2\to L^\infty}$$ $$\times\left\|\langle x\rangle^{n/2+\sigma+\epsilon}\left(\varphi(P(h))-\varphi(P_0(h))\right)\right\|_{L^1\to L^2}$$ 
$$+\left\|\langle x\rangle^{-\sigma}\sin\left(t\sqrt{P_0(h)}\right)\varphi_1^\sharp(P_0(h))\langle x\rangle^{-n/2-\sigma-\epsilon}
\right\|_{L^2\to L^\infty}$$ $$\times\left\|\langle x\rangle^{n/2+\sigma+\epsilon}\left(\varphi^\sharp(P(h))-
\varphi^\sharp(P_0(h))\right)\right\|_{L^1\to L^2}$$
$$\le C_\varepsilon h^{2\nu-n-1/2}\mu(h)^{(n+1)/2+\sigma+\varepsilon}t^{-(n-1)/2-\sigma}+Ch^{\nu-\sigma-3n/2}t^{-(n-1)/2-\sigma}.\eqno{(4.11)}$$
Furthermore, given any $f,g\in L^1$, using (2.4), (3.1) and (4.4), we get
$$t^{(n-1)/2+\sigma}\left|\left\langle\Phi_2(t,h)\langle x\rangle^{-\sigma}f,\langle x\rangle^{-\sigma}g\right\rangle\right|\le C\int_0^{t/2}(t-\tau)^{(n-1)/2+\sigma}$$ $$\times\left|\left\langle {\cal P}(h)e^{i\tau\sqrt{P(h)}}
\varphi(P(h))\langle x\rangle^{-\sigma}f,\sin\left((t-\tau)\sqrt{P_0(h)}\right)\varphi_1^\sharp(P_0(h))\langle x\rangle^{-\sigma}g\right\rangle\right|d\tau$$
$$+ C\int_{t/2}^t\tau^{(n-1)/2+\sigma}\left|\left\langle {\cal P}(h)
e^{i\tau\sqrt{P(h)}}\varphi(P(h))\langle x\rangle^{-\sigma}f,\sin\left((t-\tau)\sqrt{P_0(h)}\right)\varphi_1^\sharp(P_0(h))\langle x\rangle^{-\sigma}g\right\rangle\right|d\tau$$
 $$\le Ch^\nu\int_0^{t/2}(t-\tau)^{(n-1)/2+\sigma}\left\|\langle x\rangle^{-1/2-\epsilon}
e^{i\tau\sqrt{P_0(h)}}\varphi(P_0(h))\langle x\rangle^{-\sigma}f\right\|_{L^2}$$ 
$$\times\left\|\langle x\rangle^{-n/2-\sigma-\epsilon}\sin\left((t-\tau)\sqrt{P_0(h)}\right)\varphi_1^\sharp(P_0(h))\langle x\rangle^{-\sigma}g\right\|_{L^2}d\tau$$
 $$+ Ch^\nu\int_0^{t/2}(t-\tau)^{(n-1)/2+\sigma}\left\|\langle x\rangle^{-1-\epsilon}
\left(e^{i\tau\sqrt{P(h)}}\varphi(P(h))-e^{i\tau\sqrt{P_0(h)}}\varphi(P_0(h))\right)\langle x\rangle^{-\sigma}f\right\|_{L^2}$$ 
$$\times\left\|\langle x\rangle^{-n/2-\sigma-\epsilon}\sin\left((t-\tau)\sqrt{P_0(h)}\right)\varphi_1^\sharp(P_0(h))\langle x\rangle^{-\sigma}g\right\|_{L^2}d\tau$$
$$+ Ch^\nu\int_{t/2}^t\tau^{(n-1)/2+\sigma}\left\|\langle x\rangle^{-n/2-\sigma-\epsilon}e^{i\tau\sqrt{P_0(h)}}
\varphi(P_0(h))\langle x\rangle^{-\sigma}f\right\|_{L^2}$$ $$\times\left\|\langle x\rangle^{-1/2-\epsilon}
\sin\left((t-\tau)\sqrt{P_0(h)}\right)\varphi_1^\sharp(P_0(h))\langle x\rangle^{-\sigma}g\right\|_{L^2}d\tau$$
$$+ Ch^\nu\int_{t/2}^t\tau^{(n-1)/2+\sigma}\left\|\langle x\rangle^{-n/2-\sigma-\epsilon}\left(e^{i\tau\sqrt{P(h)}}
\varphi(P(h))-e^{i\tau\sqrt{P_0(h)}}\varphi(P_0(h))\right)\langle x\rangle^{-\sigma}f\right\|_{L^2}$$ 
$$\times\left\|\langle x\rangle^{-1-\epsilon}\sin\left((t-\tau)\sqrt{P_0(h)}\right)\varphi_1^\sharp(P_0(h))\langle x\rangle^{-\sigma}g\right\|_{L^2}d\tau$$
 $$\le Ch^\nu\left(\int_0^\infty\left\|\langle x\rangle^{-1/2-\epsilon}
e^{i\tau\sqrt{P_0(h)}}\varphi(P_0(h))\langle x\rangle^{-\sigma}f\right\|_{L^2}^2d\tau\right)^{1/2}$$ 
$$\times\left(\int_0^\infty\tau^{n-1+2\sigma}\left\|\langle x\rangle^{-n/2-\sigma-\epsilon}\sin\left(\tau\sqrt{P_0(h)}\right)
\varphi_1^\sharp(P_0(h))\langle x\rangle^{-\sigma}g\right\|_{L^2}^2d\tau\right)^{1/2}$$
 $$+ Ch^\nu\left(\int_0^\infty\tau^{n-1+2\sigma}\left\|\langle x\rangle^{-n/2-\sigma-\epsilon}
e^{i\tau\sqrt{P_0(h)}}\varphi(P_0(h))\langle x\rangle^{-\sigma}f\right\|_{L^2}^2d\tau\right)^{1/2}$$ 
$$\times\left(\int_0^\infty\left\|\langle x\rangle^{-1/2-\epsilon}\sin\left(\tau\sqrt{P_0(h)}
\right)\varphi_1^\sharp(P_0(h))\langle x\rangle^{-\sigma}g\right\|_{L^2}^2d\tau\right)^{1/2}$$
 $$+ Ch^{\nu}\left(\int_0^{t/2}\left\|\langle x\rangle^{-1-\epsilon}
\left(e^{i\tau\sqrt{P(h)}}\varphi(P(h))-e^{i\tau\sqrt{P_0(h)}}\varphi(P_0(h))\right)\langle x\rangle^{-\sigma}f\right\|_{L^2}^2d\tau\right)^{1/2}$$
 $$\times\left(\int_0^{t/2}\tau^{n-1+2\sigma}\left\|\langle x\rangle^{-n/2-\sigma-\epsilon}
\sin\left(\tau\sqrt{P_0(h)}\right)\varphi_1^\sharp(P_0(h))\langle x\rangle^{-\sigma}g\right\|_{L^2}^2d\tau\right)^{1/2}$$
 $$+ C_\varepsilon h^{2\nu-(n+1)/2}\mu(h)^{(n+1)/2+\sigma+\varepsilon}\|f\|_{L^1}\int_0^{t/2}
\left\|\langle x\rangle^{-1-\epsilon}\sin\left(\tau\sqrt{P_0(h)}\right)\varphi_1^\sharp(P_0(h))\langle x\rangle^{-\sigma}g\right\|_{L^2}d\tau$$
 $$\le \left(Ch^{\nu-\sigma-\frac{3n+1}{2}}+C_\varepsilon h^{2\nu-\sigma-\frac{3n+1}{2}}
\mu(h)^{\frac{3}{2}+\varepsilon}+ C_\varepsilon h^{2\nu-\frac{2n+3}{2}}
\mu(h)^{\frac{n+1}{2}+\sigma+\varepsilon}\right)\|f\|_{L^1}\|g\|_{L^1}.\eqno{(4.12)}$$
Clearly, (1.9) follows from (4.7), (4.11) and (4.12). The bound (1.11) can be proved in a similar way using (4.2), (4.5) and (4.10).
\eproof

{\it Proof of Theorem 1.2} Let $\varphi\in C_0^\infty((0,+\infty))$. It follows from (1.9) and (1.11) that we have the estimates
$$\left\|\langle x\rangle^{-\sigma}e^{it\sqrt{G}}\varphi(h^2G)\langle x\rangle^{-\sigma}\right\|_{L^1\to L^\infty}\le C_\varepsilon h^{-(n+1)/2-p_n(\sigma)-\varepsilon/2}|t|^{-(n-1)/2-\sigma},
\eqno{(4.13)}$$
$$\left\|\langle x\rangle^{-\sigma}e^{itG}\varphi(h^2G)\langle x\rangle^{-\sigma}\right\|_{L^1\to L^\infty}\le C_\varepsilon h^{\sigma-q_n(\sigma)-\varepsilon/2}|t|^{-n/2-\sigma}.\eqno{(4.14)}$$
We now write
$$\left(\sqrt{G}\right)^{-(n+1)/2-p_n(\sigma)-\varepsilon}\chi\left(\sqrt{G}\right)=\int_0^1\psi\left(h\sqrt{G}\right)
h^{(n+1)/2+p_n(\sigma)+\varepsilon-1}dh,\eqno{(4.15)}$$
where $\psi(\lambda)=\lambda^{1-(n+1)/2-p_n(\sigma)-\varepsilon}\chi'(\lambda)\in C_0^\infty((0,+\infty))$. 
By (4.13) and (4.15) we get
$$\left\|\langle x\rangle^{-\sigma}e^{it\sqrt{G}}\left(\sqrt{G}\right)^{-(n+1)/2-p_n(\sigma)-\varepsilon}\chi\left(\sqrt{G}\right)\langle x\rangle^{-\sigma}\right\|_{L^1\to L^\infty}$$ 
$$\le\int_0^1\left\|\langle x\rangle^{-\sigma}e^{it\sqrt{G}}\psi\left(h\sqrt{G}\right)\langle x\rangle^{-\sigma}\right\|_{L^1\to L^\infty}h^{(n+1)/2+p_n(\sigma)+\varepsilon-1}dh$$ $$
\le C_\varepsilon |t|^{-(n-1)/2-\sigma}\int_0^1h^{-1+\varepsilon/2}dh\le C'_\varepsilon |t|^{-(n-1)/2-\sigma}.$$
Similarly, using (4.14) we get
$$\left\|\langle x\rangle^{-\sigma}e^{itG}\left(\sqrt{G}\right)^{\sigma-q_n(\sigma)-\varepsilon}\chi\left(\sqrt{G}\right)\langle x\rangle^{-\sigma}\right\|_{L^1\to L^\infty}$$ 
$$\le\int_0^1\left\|\langle x\rangle^{-\sigma}e^{itG}\psi_1\left(h\sqrt{G}\right)\langle x\rangle^{-\sigma}\right\|_{L^1\to L^\infty}h^{-\sigma+q_n(\sigma)+\varepsilon-1}dh$$ $$
\le C_\varepsilon |t|^{-n/2-\sigma}\int_0^1h^{-1+\varepsilon/2}dh\le C'_\varepsilon |t|^{-n/2-\sigma}.$$
To prove (1.18) observe that, if $k<1$ and $\delta$ and $\sigma$ satisfy the conditions of Theorem 1.2, we have $\sigma>q_n(\sigma)$. Hence
$$\left\|\langle x\rangle^{-\sigma}e^{itG}\chi\left(\sqrt{G}\right)\langle x\rangle^{-\sigma}\right\|_{L^1\to L^\infty}$$ 
$$\le\int_0^1\left\|\langle x\rangle^{-\sigma}e^{itG}\psi_2\left(h\sqrt{G}\right)\langle x\rangle^{-\sigma}\right\|_{L^1\to L^\infty}h^{-1}dh$$ $$
\le C_\varepsilon |t|^{-n/2-\sigma}\int_0^1h^{\sigma-q_n(\sigma)-1-\varepsilon/2}dh\le C'_\varepsilon |t|^{-n/2-\sigma}.$$
\eproof

\section*{Appendix}

In this appendix we will sketch the proof of the estimates (4.3)-(4.6). To this end, we will use the fact that the kernels of the operators 
$e^{it\sqrt{P_0(h)}}\varphi(P_0(h))$ and $e^{itP_0(h)}\varphi(P_0(h))$ are of the form $K_h(|x-y|,t)$ and $\widetilde K_h(|x-y|,t)$, respectively, where
$$K_h(w,t)=\frac{w^{2-n}}{(2\pi)^{n/2}}\int_0^\infty e^{ith\lambda}\varphi(h^2\lambda^2){\cal J}_{\frac{n-2}{2}}(w\lambda)\lambda d\lambda=h^{-n}K_1(w/h,t),\eqno{(A.1)}$$
$$\widetilde K_h(w,t)=\frac{w^{2-n}}{(2\pi)^{n/2}}\int_0^\infty e^{ith^2\lambda^2}\varphi(h^2\lambda^2){\cal J}_{\frac{n-2}{2}}(w\lambda)\lambda d\lambda=h^{-n}\widetilde K_1(w/h,t),\eqno{(A.2)}$$
where ${\cal J}_{\frac{n-2}{2}}(z)=z^{(n-2)/2}J_{\frac{n-2}{2}}(z)$, $J_{\frac{n-2}{2}}(z)$ being the Bessel function of order $\frac{n-2}{2}$.
In view of the inequality
$$\langle x\rangle^{-\sigma}\langle y\rangle^{-\sigma}\le \langle x-y\rangle^{-\sigma},\quad\forall\sigma\ge 0,$$
it is easy to see that the estimates (4.3)-(4.6) follow from the following

\quad

\noindent
{\bf Lemma A.1.} {\it For all $w>0, t\neq 0, 0<h\le 1, s\ge 0$, we have}
$$\left|K_h(w,t)\right|\le C|t|^{-s}h^{-s-(n+1)/2}g_s(w),\eqno{(A.3)}$$
$$\int_{-\infty}^\infty|t|^{2s}\left|K_h(w,t)\right|^2dt\le Ch^{-2s-n-1}g_s(w)^2,\eqno{(A.4)}$$
$$\left|\widetilde K_h(w,t)\right|\le C|t|^{-s-1/2}h^{-s-(n+1)/2}g_s(w),\eqno{(A.5)}$$
$$\int_{-\infty}^\infty|t|^{2s}\left|\widetilde K_h(w,t)\right|^2dt\le Ch^{-2s-n-1}g_s(w)^2,\eqno{(A.6)}$$
{\it where $g_s(w)=w^{s-(n-1)/2}$ if $s\le (n-1)/2$, $g_s(w)=\langle w\rangle^{s-(n-1)/2}$ if $s\ge (n-1)/2$.}

\quad

{\it Proof.} In view of the identities (A.1) and (A.2), it is clear that it suffices to prove (A.3)-(A.6) for $h=1$. Let first $w\le 1$.
Recall that near $z=0$ the function ${\cal J}_{\frac{n-2}{2}}(z)$ is equal to $z^{n-2}$ times an analytic function. Using this and integrating by parts, it is easy to see that in this case the functions $K_1$ and $\widetilde K_1$ satisfy the bounds
$$\left|K_1(w,t)\right|\le C|t|^{-s},\eqno{(A.7)}$$
$$\left|\widetilde K_1(w,t)\right|\le C|t|^{-s-1/2},\eqno{(A.8)}$$
for every $s\ge 0$. Clearly, when $w\le 1$ the estimates (A.3)-(A.6) follow from (A.7) and (A.8). Let now $w\ge 1$. In this case we will use the fact that for $z\gg 1$ the function ${\cal J}_{\frac{n-2}{2}}(z)$ is of the form $e^{iz}b^+(z)+e^{-iz}b^-(z)$, where $b^\pm(z)$ are symbols of order $\frac{n-3}{2}$. Given any integers $k,\ell\ge 0$, set
$$b_k^\pm(z)=e^{\mp iz}\frac{d^k}{dz^k}\left(e^{\pm iz}b^\pm(z)\right),\quad b_{k,\ell}^\pm(z)=\frac{d^\ell}{dz^\ell}b_k^\pm(z).$$
Clearly, $b_k^\pm(z)$ are also symbols of order $\frac{n-3}{2}$. Hence
$$\left|b_{k,\ell}^\pm(z)\right|\le C_{k,\ell}z^{\frac{n-3}{2}-\ell},\quad\forall z\ge 1.\eqno{(A.9)}$$
Let $m,N\ge 0$ be integers. Integrating $m$ times by parts, we can write
$$K_1(w,t)=\frac{w^{2-n}}{(2\pi)^{n/2}}(it)^{-m}\sum_\pm\sum_{k=0}^m\int_0^\infty e^{i(t\pm w)\lambda}w^kb_k^\pm(w\lambda)\varphi_{k,m}(\lambda)d\lambda,$$
with some functions $\varphi_{k,m}\in C_0^\infty((0,+\infty))$ independent of $w$ and $t$. We now integrate $N$ times by parts to obtain
$$K_1(w,t)=\frac{w^{2-n}}{(2\pi)^{n/2}}(it)^{-m}\sum_\pm\sum_{k=0}^m\sum_{\ell =0}^N(t\pm w)^{-N}\int_0^\infty e^{i(t\pm w)\lambda}w^kb_{k,\ell}^\pm(w\lambda)\varphi_{k,\ell,m,N}(\lambda)d\lambda.$$
Hence, in view of (A.9), we get the bound
$$\left|K_1(w,t)\right|\le C_{m,N}w^{m-\frac{n-1}{2}}|t|^{-m}\left(|t-w|^{-N}+|t+w|^{-N}\right).\eqno{(A.10)}$$
By interpolation, (A.10) holds for all real $m\ge 0$. It is easy to see now that the estimates (A.3) and (A.4) (with $h=1$) follow from (A.10).

Integrating by parts $m$ times with respect to the variable $\lambda^2$ we can write the function $\widetilde K_1$ as follows
$$\widetilde K_1(w,t)=\frac{w^{2-n}}{(2\pi)^{n/2}}(it)^{-m}\sum_{k=0}^m\int_0^\infty e^{it\lambda^2}\widetilde\varphi_{k,m}(\lambda)\frac{d^k}{d(\lambda^2)^k}{\cal J}_{\frac{n-2}{2}}(w\lambda)d\lambda$$
$$=\frac{w^{2-n}}{(2\pi)^{n/2}}(it)^{-m}\sum_{k=0}^m\int_0^\infty e^{it\lambda^2}f_{k,m}(w,\lambda)d\lambda,$$
where
$$f_{k,m}(w,\lambda)=\widetilde\varphi_{k,m}^\sharp(\lambda)\frac{d^k}{d\lambda^k}{\cal J}_{\frac{n-2}{2}}(w\lambda).$$
We now apply the inequality
$$\left|\int_0^\infty e^{it\lambda^2}f(\lambda)d\lambda\right|\le C|t|^{-1/2}\left\|\widehat f\right\|_{L^1},\quad\forall f\in C_0^\infty({\bf R}),$$
to get
$$\left|\widetilde K_1(w,t)\right|\le C_mw^{2-n}|t|^{-m-1/2}\sum_{k=0}^m\left\|\widehat f_{k,m}(\cdot,w)\right\|_{L^1}.\eqno{(A.11)}$$
On the other hand, as above one can see that the function $\widehat f_{k,m}$ satisfies the bound
$$\left|\widehat f_{k,m}(\tau,w)\right|\le C_{N}w^{k+\frac{n-3}{2}}\left(|\tau-w|^{-N}+|\tau+w|^{-N}\right)\eqno{(A.12)}$$
for every integer $N\ge 0$. By (A.12)
$$\left\|\widehat f_{k,m}(\cdot,w)\right\|_{L^1}\le Cw^{k+\frac{n-3}{2}}.\eqno{(A.13)}$$
By (A.11) and (A.13)
$$\left|\widetilde K_1(w,t)\right|\le C_mw^{m-\frac{n-1}{2}}|t|^{-m-1/2}\eqno{(A.14)}$$
for every integer $m\ge 0$, and hence by interpolation for all real $m\ge 0$, which in turn proves (A.5). 
It is easy also to see that (A.6) (with $h=1$) follows from (A.14). Indeed, applying (A.14) with $m=s-\epsilon$ and $m=s+\epsilon$, we have
$$\int_{-\infty}^\infty|t|^{2s}\left|\widetilde K_h(w,t)\right|^2dt$$ $$\le Cw^{2s-2\epsilon-n+1}\int_{|t|\le w}|t|^{-1+2\epsilon}dt+
Cw^{2s+2\epsilon-n+1}\int_{|t|\ge w}|t|^{-1-2\epsilon}dt\le Cw^{2s-n+1}.$$

\eproof

{\bf Acknowledgements.} A part of this work has been carried out while G. V. was visiting the Universidade Federal de Pernambuco, 
Brazil, in January, February 2011 with the partial support of the agreement Brazil-France in Mathematics - Proc. 69.0014/01-5.
The first two authors have been partially supported by the CNPq-Brazil.

F. Cardoso

Universidade Federal de Pernambuco, 

Departamento de Matem\'atica, 

CEP. 50540-740 Recife-Pe, Brazil,

e-mail: fernando@dmat.ufpe.br

\quad

C. Cuevas

Universidade Federal de Pernambuco, 

Departamento de Matem\'atica, 

CEP. 50540-740 Recife-Pe, Brazil,

e-mail: cch@dmat.ufpe.br

\quad

G. Vodev

Universit\'e de Nantes,

 D\'epartement de Math\'ematiques, UMR 6629 du CNRS,
 
 2, rue de la Houssini\`ere, BP 92208, 
 
 44332 Nantes Cedex 03, France,
 
 e-mail: vodev@math.univ-nantes.fr

\end{document}